\documentclass[a4paper,10pt]{article}

\usepackage{mathrsfs}

\usepackage{amssymb}
\usepackage{hyperref} 
\date{}
\input xy 
\xyoption{all}  
\input{xypic} 

\renewcommand{\geq}{\geqslant}

\newcommand{\sch}[1]{\mathscr #1}

\newcommand{\etal}{_{\rm \acute{e}t}}

\newcommand{\hotimes}{\widehat{\otimes}}
\renewcommand{\phi}{\varphi}
\renewcommand{\dim}[1]{\mathsf{dim}_{#1}\;}

\renewcommand{\Bbb}{\mathbb}

\newcommand{\shl}[1]{\mathsf{Sh}(#1,\ell)}

\newcommand{\R}{\mbox{\bf R}}
\newcommand{\gm}{{\Bbb G}_{m}}

\renewcommand{\H}{\mbox{\rm H}}
\renewcommand{\R}{{\rm R}}
\newcommand{\RR}{{\Bbb R}}
\newcommand{\HH}{{\Bbb H}}

\newcommand{\ZZ}{{\Bbb Z}}

\newcommand{\FF}{{\Bbb F}}
\newcommand{\DD}{{\Bbb D}}
\newcommand{\CC}{{\Bbb C}}
\newcommand{\PP}{{\Bbb P}}
\newcommand{\Aff}{{\Bbb A}}

\newcommand{\NN}{{\Bbb N}}

\newcommand{\QQ}{{\Bbb Q}}

\renewcommand{\epsilon}{\varepsilon}

\newcommand{\an}{^{an}}

\newcommand{\spec}{\mathsf{Spec}\;}

\newcommand{\red}{\widetilde}

\newcommand{\hres}{{\sch H}}

 \title{\'Etale cohomology of schemes and analytic spaces}
 \author{Antoine Ducros}
 \begin{document}
 \maketitle
 \setcounter{section}{-1}
 
 \section{Introduction}
 
  This text corresponds essentially to the four hour mini-course that was given by the author in July 2010 at the university Paris 6 (Jussieu) during the {\em Summer school on Berkovich spaces} which was organized by the ANR Project {\em Espaces de Berkovich}. 
 
 The goal of those talks was to make a survey on scheme-theoretic and Berkovich-analytic \'etale cohomology theories for people coming from various mathematical areas (including model theory, complex dynamics, tropical geometry...) and not necessarily familiar with algebraic geometry {\em \`a la} Grothendieck. That's the reason why we chose to rather insist on the motivation, the rough ideas behind the theory and the basic examples, than to give proofs -- those are usually technically involved, and can be found in the literature. Nevertheless, the interested reader may wish to have some technical references about this subject; let me now give some of them.

\medskip
\noindent
{\bf About the notion of a derived functor}. The theory was written down for the first time by Grothendieck in a very beautiful, and quite easily readable, paper (\cite{thk}). 
 
 \medskip
 \noindent
 {\bf About scheme-theoretic \'etale cohomology}. Foundations are developped in \cite{sga1}, \cite{sga41}-\cite{sga43}, and \cite{sga5}. More precisely, \cite{sga1} introduces the notion of an \'etale morphism, and Grothendieck's theory of the fundamental group; the three books \cite{sga41}-\cite{sga43} are devoted to Grothendieck topologies, general properties of the associated cohomology theories, and the \'etale cohomology of schemes; and the $\ell$-adic formalism is introduced in \cite{sga5}. 

\medskip
But all those books, though extremely interesting and useful, are quite difficult to read because of their high degree of generality. Thus the reader may prefer have a look to \cite{sga45}, in which the theory is presented without (too much detailed) proofs, and with a lot of explanations about the ideas, the geometric interpretation, the analogy with usual topology and so on; he could as well read \cite{mil}, which is certainly far more understandable than the SGA's, though it also requires some familiarity with schemes -- but such a familiarity is however unavoidable for anybody who wants to work seriously with \'etale cohomology. 

\medskip
Let us also mention a beautiful text by Luc Illusie (\cite{ill}) about the \'etale cohomology from a double viewpoint, historical and mathematical. 

\medskip
\noindent
{\bf About analytic \'etale cohomology.} The fundamental reference for that subject is the foundational article by Berkovich himself (\cite{brk2}), which is not easily readable by somebody who doesn't already master basic analytic geometry and Grothendiecks' language (Grothendieck topologies, sheaves, derived categories and so on). The comparison theorem of \cite{brk2} has been improved in \cite{brkc}, which is also quite technical. 

\medskip
\noindent
{\bf Acknowledgements.} I would like to thank Mathilde Herblot for having written detailed notes during the lectures, and allowed me to use them for the redaction of this text. 
  
 \section{The general motivation} 
 
 As well for schemes as for Berkovich analytic spaces, \'etale cohomology is a cohomology theory which behaves like singular cohomology of complex analytic spaces; that is to say, it does satisfy the same general theorems (bounds on cohomological dimension,  K\"unneth formula, Poincar\'e duality...), and 'takes the expected values' on several classes of algebraic varieties (projective curves, abelian varieties) or analytic spaces (polydiscs, polyannuli, analytic tori....). Moreover, as we will see, this theory also encodes in a very natural way the Galois theory of fields, which helps thinking to the latter, which is {\em a priori} purely algebraic, in a quite geometrical way.  
 
 \medskip
In both scheme-theoretic and analytic contexts, the need for a 'singular-like' cohomology theory came from very deep {\em arithmetical} conjectures:

\medskip
- as far as schemes are concerned, the Weil conjectures; Weil himself suggested to prove them by defining such a cohomology, which is for that reason often called a {\em Weil cohomology}, and by applying it to algebraic varieties defined over finite fields; it is for that purpose that Grothendieck {\em et. al.} developped the theory in the sixties, and Weil conjectures were eventually proven by Deligne in 1974; 

-as far as analytic spaces are concerned, the Langlands program; Drinfeld and Carayol had conjectured that the local Langlands correspondence should 'live' in the \'etale cohomology vector spaces of some meaningful $p$-adic analytic moduli spaces, whose archetype is the so-called {\em $p$-adic upper half plane} $\PP^{1,an}_{\QQ_p}-\PP^1(\QQ_p)$;  but at that time, the \'etale cohomology of a $p$-adic analytic space was just not defined. It was done for the first time by Berkovich, who used it to prove a conjecture by Deligne on vanishing cycles; after that, it was shown to have the expected applications to the local Langlands program (works by Harris, Taylor, Boyer, Hausberger, Dat....). 

\medskip
{\small 
\noindent
{\bf The $p$-adic upper half-plane}. Remind that Berkovich analytic spaces are defined over {\em any} (i.e., possibly Archimedean) complete valued field. In particular, $\Aff^{1,an}_{\RR}$ does make sense; as a topological space, it is the set of multiplicative semi-norms on $\RR[T]$ which extend the usual absolute value on $\RR$. The map $z\mapsto (P\mapsto |P(z)|)$ is easily seen to induce an homeomorphism, and even an isomorphism of locally ringed spaces, between  the quotient of $\CC$ by the complex conjugation and $\Aff^{1,an}_{\RR}$.

The space $\PP^{1,an}_{\RR}$ is obtained by glueing two copies of $\Aff^{1,an}_{\RR}$ in the usual way; hence there is a canonical isomorphism between $\PP^{1,an}_{\RR}$ and the quotient of $\PP^1(\CC)$ by the complex conjugation; therefore, $(\PP^{1,an}_{\RR}-\PP^1(\RR))\simeq \HH$, where $\HH$ denotes the usual upper half-plane. 

\medskip
This is one of the numerous reasons why $\PP^{1,an}_{\QQ_p}-\PP^1(\QQ_p)$ should be thought as the $p$-adic analog of $\HH$. Note that it is a well-defined Berkovich space: indeed, as $\QQ_p$ is locally compact, $\PP^1(\QQ_p)$ is compact; as a consequence, $\PP^{1,an}_{\QQ_p}-\PP^1(\QQ_p)$ is an open subset of $\PP^(\QQ_p)$ (one checks that it is dense, and connected); thus it inherits a natural structure of a $p$-adic analytic space. }

\medskip
Let us make now two remarks. 

\medskip
\noindent
{\bf First remark.} In both scheme-theoretic and analytic settings, there are natural cohomology theories, namely the ones associated with the {\em topological spaces} involved. But though they are highly interesting, they definitely can {\em not} play the expected role. Let us give one reason for that : if $X$ is a $d$-dimensional variety (resp. a $d$-dimensional paracompact analytic space) over an algebraically closed field, then its topological cohomology vanishes in rank $\geq d+1$, while non-zero cohomology groups should be allowed to occur up to rank $2d$ in a singular-like cohomology theory. 

\medskip
\noindent
{\bf Second remark (due to Serre)}. Let us explain why it is not possible to hope for a singular-like cohomology theory for schemes providing $\QQ$-vector spaces (the same argument would also work {\em mutatis mutandis} for analytic spaces). Assume that there exists such a theory $\H^\bullet$, let us fix a prime number $p$ and let $\overline {\FF_p}$ be an algebraic closure of $\FF_p$. There exists an elliptic curve $X$ over $\overline{\FF_p}$ with the following property: {\em the ring $(\mathsf {End}\;X)\otimes_{\ZZ}\QQ$ is a non-split $\QQ$-quaternion algebra, that is, a non-commutative field with center $\QQ$ and whose dimension over $\QQ$ is equal to 4}\footnote{It can be shown that $X$ satisfies this property if and only if it is {\em supersingular}, that is, if and only if its $p$-torsion subgroup consists in {\em one} $\overline{\FF_p}$-point, with multiplicity $p^2$.}. By our assumptions on $\H^\bullet$, the group $\H^1(X)$ is a $2$-dimensional $\QQ$-vector space, on which $\mathsf {End} \;X$ acts. The corresponding homomorphism $$\mathsf {End} \;X\to  \mathsf{End} \;\H^1(X)\simeq \mathsf M_2(\QQ)$$ induces an arrow $$(\mathsf {End} \;X)\otimes_{\ZZ}\QQ\to  \ \mathsf M_2(\QQ).$$ As the left-hand side ring is a (non-commutative) field, this map is injective. By comparison of dimensions it is bijective; but $\mathsf M_2(\QQ)$ is not a field, contradiction. 
 
 \medskip
 This is the reason why \'etale cohomology usually deals with torsion or $\ell$-adic coefficients, and (almost) never with integral or rational ones.  
 
\section{The notion of a Grothendieck topology}\label{grotop}

{\bf \em General references: {\rm  \cite{sga41} } or {\rm \cite{mil}}, far easier to read.} 
 
\medskip
Classical general topology consists in axiomatizing the notion of an open subset (or of a neighborhood); in order to develop \'etale cohomology, Grothendieck slightly generalized this formalism, extending it to the categorical setting, and giving rise to what one calls a {\em Grothendieck topology}. In that approach,  one 'only' axiomatizes the notion of a {\em covering} -- this is precisely what is needed to define sheaves and their cohomology groups. 
 
 \medskip
 Now let us start with a category $\mathsf C$. We will note give here the precise definition of a category; let us simply remind that one has {\em objects} of $\mathsf C$, {\em morphisms} between them (including the identity map of every object) and that there is a notion of composition of morphisms which satisfies the usual rules. 
 
 \medskip
 We will moreover assume that if $Y\to X$ and $Z\to X$ are arrows of $\mathsf C$, the {\em fiber product} $Y\times_XZ$ does exist in $\mathsf C$; let us give some explanations. One says that an object $T$ of $\mathsf C$, equipped with two morphisms $T\to Y$ and $T\to Z$ uch that the diagram $\diagram T\rto\dto&Z\dto\\Y\rto& X\enddiagram$ commutes, is a {\em fiber product of $Y$ and $Z$ over $X$} if for every object $S$ of $\mathsf C$, and every couple of morphisms $(S\to Y,S\to Z)$ such that $\diagram S\rto\dto&Z\dto\\Y\rto& X\enddiagram$ commutes, there exists a unique morphism $S\to T$ making the diagram $$
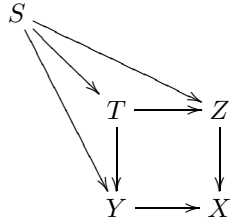
\diagram S\drto \drrto \ddrto& & \\ & T\rto\dto&Z\dto\\&Y\rto& X\enddiagram$$ commute. 
 
{\em As soon as such a fiber product exists}, it is unique up to canonical isomorphism; it is therefore called {\em the} fiber product of $Y$ and $Z$ over $X$, and denoted by $Y\times_XZ$. 

\medskip
{\bf Examples of categories in which $Y\times_XZ$ always exists.}

\medskip
$\bullet$ {\em The category of sets.} In this case $Y\times_XZ$ is the set-theoretic fiber product, that is, the set of couples $(y,z)$ with $y\in Y$ and $z\in Z$ such that $y$ and $z$ have the same image on $X$. 

$\bullet$ {\em The category of topological spaces.} In this case $Y\times_XZ$ is the set-theoretic fiber product {\em endowed with the product topology}. 

$\bullet$ {\em The category of schemes.} Let us assume for simplicity that $X,Y$ and $Z$ are affine, say respectively equal to $\spec A, \spec B$ and $\spec C$. Then $Y\times_XZ$ is nothing but $\spec (B\otimes_AC)$. 

$\bullet$ {\em The category of analytic spaces over a non-Archimedean, complete field.} Let us assume for simplicity that $X,Y$ and $Z$ are affinoid, say respectively equal to $\sch M (A), \sch M(B)$ and $\sch M( C)$. Then $Y\times_XZ$ is nothing but $\sch M (B\hotimes_AC)$. 

$\bullet$ {\em The category of the open subsets of a given topological space, in which morphisms are the inclusions.} In that case, $Y\times_XZ$ is nothing but the intersection $Y\cap Z$.

$\bullet$ {\em The category of analytic domains of a given analytic space, in which morphisms are the inclusions.} In that case, $Y\times_XZ$ is nothing but the intersection $Y\cap Z$.

\medskip
Now, let us go back to our category $\mathsf C$ in which the fiber products are assumed to exist. A {\em Grothendieck topology} on $\mathsf C$ consists in the datum, for every object $X$ of $\mathsf C$, of a class of families of arrow $(X_i\to X)_i$, which are called {\em coverings} of $X$, and which are subject to the following axioms: 

\medskip
1) If $Y\to X$ is an isomorphism, it is a covering of $X$; 

2) if $(X_i\to X)_i$ is a covering of $X$, and if $Y\to X$ is any morphism, $$(X_i\times_XY)_i\to Y$$ is a covering of $Y$; 

3) if $(X_i\to X)_i$ is a covering of $X$ and if $(X_{ij}\to X_i)_j$ is for every $i$ a covering of $X_i$ then $(X_{ij}\to X)_{i,j}$ is a covering of $X$. 

\medskip
For the moment, we will only give two examples, which are quite set-theoretic; we will define further the \'etale topology on a scheme (or an analytic space) which will be far more categorical. 

\medskip
\noindent
{\bf First example.} On the category of open subsets of a given topological space, the usual open coverings define a Grothendieck topology. 

\medskip
\noindent
{\bf Second example.} On the category of analytic domains of a given analytic space, the admissible coverings define a Grothendieck topology; it is called the {\em G-topology.}

\medskip
\noindent
{\bf Definition.} A {\em site} is a category equipped with a Grothendieck topology. 

\subsection*{Presheaves and sheaves} 

We fix from now on a site $\mathsf C$. 

\medskip
\noindent
{\bf Comments on the vocabulary.} Though the behavior of $\mathsf C$ can be slightly different from that of an honest topological space, one will often use, while working with $\mathsf C$, the vocabulary and, to some reasonable extent, the intuition of classical topology. Let us give a quite vague example: if $X$ is an objet of $\mathsf C$ and if $\mathsf P$ is a property, we will say that $\mathsf P$ is true {\em locally on $X$} if there exists a covering $(X_i\to X)$ such that $\mathsf P$ is true on every $X_i$; we will say that $\mathsf P$ is local if it is true on a given object of $\mathsf C$ as soon as it is true locally on it. 

\medskip
\noindent
{\bf The notion of a presheaf}. A {\em presheaf} $\sch F$ on $\mathsf C$ is a {\em contravariant functor} from $\mathsf C$ to the category of sets; that is, an assignment which sends any object $X$ of $\mathsf C$ to a set $\sch F(X)$, together with the datum, for every arrow $Y\to X$ in $\mathsf C$, of a {\em restriction map} $\sch F(X)\to \sch F(Y)$ (note that it goes 'the wrong way'), all of this having to behave well with respect to the composition of arrows. The elements of $\sch F(X)$ for a given $X$ will sometimes be called {\em sections of $\sch F$ on $X$.}

If there is no risk of confusion\footnote{There could be one because $\sch F(X)\to \sch F(Y)$ not only depends on $Y$, but on the arrow $Y\to X$.}, the restriction $\sch F(X)\to \sch F(Y)$ will be denoted $s\mapsto s_{|Y}$.

\medskip
A {\em map} from a presheaf $\sch F$ to a presheaf $\sch G$ consists in a family of maps $\sch F(X)\to\sch G(X)$, where $X$ goes through the collection of all objects of $\mathsf C$, which is compatible with restriction maps.

\medskip
\noindent
{\bf The notion of a sheaf}. A presheaf $\sch F$ on $\mathsf C$ is said to be a {\em sheaf} if it satisfies the following condition: {\em for every object $X$ of $\mathsf C$, for every covering $(X_i\to X)$, and for every family $(s_i)\in \prod\limits_i \sch F(X_i)$ such that $$s_{i|X_i\times_XX_j}=s_{j|X_i\times_XX_j}\in \sch F(X_i\times_X X_j)$$ for all $(i,j)$, there exists a unique $s\in \sch F(X)$ such that $s_{|X_i}=s_i$ for every $i$.}

\medskip
The section $s$ shoulds be thought as beeing obtained by {\em glueing} the $s_i$'s, the equalities $s_{i|X_i\times_XX_j}=s_{j|X_i\times_XX_j}$ beeing the coincidence conditions that make the glueing process possible.

\medskip
If $\sch F$ and $\sch G$ are two sheaves on $\mathsf C$,  a map from $\sch F$ to $\sch G$ is simply a map between the {\em presheaves} $\sch F$ and $\sch G$. 

\medskip
\noindent
{\bf Link with old-fashioned sheaves.} If $\mathsf C$ is the category of all open subsets of a given topological space $T$, equipped with the aforementioned Grothendieck topology, then a sheaf on $\mathsf C$ is nothing but a classical sheaf on the topological space $T$. 

\medskip
\noindent
{\bf Presehaves and sheaves with some extra-structures.} We have actually just defined presheaves and sheaves of {\em sets.} We will also use the notion of presheaves and sheaves of abelian groups. The definitions are almost the same, except that $\sch F(X)$ comes for every object $X$ of $\mathsf C$ with the structure of an abelian group and that restrictions map are required to be homomorphisms of groups ; a morphism between two presheaves or sheaves of abelian groups $\sch F$ and $\sch G$ is a map from $\sch F$ to $\sch G$ such that $\sch F(X)\to \sch G(X)$ is a morphism of groups for every object $X$ of $\mathsf C$. 

One can of course define in an analogous way presheaves and sheaves of rings, of $k$-vector spaces or $k$-algebras (for $k$ a field), etc. , and morphisms between them. 

\medskip
\noindent
{\bf Subsheaf of a sheaf.} If $\sch F$ is a sheaf on $\mathsf C$, a {\em subsheaf} of $\sch F$ is a sheaf $\sch G$ such that $\sch G(X)$ is for every object $X$ of $\mathsf C$ a subset of $\sch F(X)$. If $\sch G$ is a subsheaf of $\sch F$, the property for a section of $\sch F$ to be a section of $\sch G$ is a local property. 

\medskip
\noindent
{\bf Presheaf associated with a sheaf}. Every presheaf $\sch F$ on $\mathsf C$ can be sheafified in a natural way;  if $A$ is a set and if $\sch F$ is the constant presheaf $U\mapsto A$, its sheafification will be denoted by $\underline A_{\mathsf C}$ (or $\underline A$, or $A_C$, or $A$ if there is no risk of confusion), and called the {\em constant sheaf} (associated with) $A$. When $\mathsf C$ is the site of open subsets of a given topological space $X$, the sheaf $\underline A$ is the sheaf of locally constant $A$-valued functions on $X$. 

\medskip
If $\sch O$ is a sheaf of commutative rings on $\mathsf C$ and if $\sch F$ and $\sch G$ are two sheaves of $\sch O$-modules, the presheaf $U\mapsto \sch F(U)\otimes_{\sch O(U)}\sch G(U)$ is not a sheaf in general; its sheafification is denoted by $\sch F\otimes_{\sch O}\sch G$. 

\medskip
\noindent
{\bf A crucial notion: the {\em sheaf theoretic} image of a map.} Let $\sch F$ and $\sch G$ be two sheaves on $\mathsf C$, and let $f$ be a map from $\sch F$ to $\sch G$; one would like to define a subsheaf of $\sch G$ which could be reasonably thought as beeing the image of $f$. The problem is that the 'naive' image of $f$, that is, the assignment $$X\mapsto\{t\in \sch G(X), \exists \;s\in \sch F(X)\;{\rm s.t.}\; t=f(s)\}, $$ is a presheaf {\em which has no reason to be a sheaf}. The point is the following: for a section of $\sch G$, belonging to the naive image of $f$, is {\em not} local in general (an example is given a little bit further). 

\medskip
To bypass this difficulty, one defines the {\em sheaf-theoretic image} of the map $f$ as the assignment that sends an object $X$ of $\mathsf C$ to the set of sections of $\sch G$ on $X$ which belong {\em locally} to the naive image of $f$; in more precise words, $X$ is sent to the set of $t\in \sch G(X)$ satisfying the following properties: {\em there exists a covering $(X_i)$ of $X$ and, for every $i$, an element $s_i$ of $\sch F(X_i)$ such that $f(s_i)=t_{|X_i}$ for all $i$.} 

This assignment turns out to be a sheaf we will denote by $\mathsf{Im}\; f$; it can be shown that this is the sheafification of the naive image of $f$; this is also the smallest subsheaf of $\sch G$ containing it. We will say that $f$ is {\em surjective} if $\mathsf{Im}\;f=\sch G$. 

\medskip
\noindent
{\bf A fundamental example.} Let $\sch O$ be the sheaf of holomorphic functions on $\CC$. The fact for a holomorphic function to admit a primitive on a given open subset of $\CC$ is not a local property -- for instance, $1/z$ has locally a primitive on $\CC^*$, but not globally; hence the naive image of the derivation $d: \sch O\to \sch O$ is not a subsheaf of $\sch O$.

However, the sheaf-theroetic image $\mathsf {Im}\;d$ of $d$ is very easy to describe: as every holomorphic function on an open subset of $\CC$ has locally a primitive, $\mathsf {Im}\;d$ is nothing but the whole sheaf $\sch O$. 

\medskip
\noindent
{\bf Exact sequences of sheaves.} Now let $\sch F,\sch F'$ and $\sch F''$ be three sheaves of abelian groups on $\mathsf C$, and let $f: \sch F'\to \sch F$ and $g: \sch F\to \sch F''$ be two morphisms. We will say that the sequence $$\diagram 0 \rto & \sch F'\rto^f&\sch F\rto ^g&\sch F''\rto & 0\enddiagram$$ is exact if the following conditions hold:

\medskip
i) for every object $X$ of $\mathsf C$, the sequence $0\to \sch F'(X)\to \sch F(X)\to \sch F''(X)$ is exact (in the usual sense for a sequence of abelian groups); 

ii) $\sch F''=\mathsf {Im}\; g$ (in other words: $g$ is surjective). 

\medskip
\noindent
{\bf Comment.} Note that in i) we {\em do not} require $\sch F(X)\to \sch F''(X)$ to be surjective for every $X$; we simply demand in ii) that $g$ is surjective, which is slightly weaker. 

\medskip
\noindent
{\bf Example.} Let $\sch O$ be the sheaf of holomorphic functions on $\CC$ and let $d$ be the derivation; remind that $\underline{\CC}$ denotes the constant sheaf associated with $\CC$, that is, the sheaf of locally constant complex-valued functions, which embeds naturally in $\sch O$. One immediatly checks that the sequence $$\diagram 0 \rto & \underline{\CC} \rto&\sch O\rto ^d&\sch O\rto & 0\enddiagram$$ is exact (the fact that $\mathsf {Im}\;d=\sch O$ has be mentioned above); nevertheless, note that the derivation $\sch O(\CC^*)\to \sch O(\CC^*)$ is {\em not} onto: again, consider the function $1/z$. 

\section{Cohomology of sheaves}\label{coho}

{\bf \em General references: {\rm \cite{sga42}} and {\rm \cite{mil}}, and {\rm \cite{thk}} for the notion of a derived functor.}

\subsection*{Abelian categories, left-exact functors, and derived functors}

We will not give here the precise definition of an {\em abelian category}. Let us just say that this is a category $\mathsf A$ which 'looks like' the category of modules over a ring; as an example, in such a category:

\medskip
$\bullet$ the set of morphisms between two given objects is an abelian group ;

$\bullet$ there is a zero object; 

$\bullet$ one can define the direct sum of a family of objects; 

$\bullet$ one can define the sub-objects of a given object, and the quotient of an object by one of its sub-objects; 

$\bullet$ one can define the image, the kernel, and the cokernel of a map, and hence the notions of surjective and an injective maps; 

$\bullet$ one can define exact sequences. 

\medskip
Moreover, all those notions behave the usual way. 

\medskip
\noindent
{\bf Examples.} Let us give two typical examples of abelian categories: the category of modules over a given ring (non necessarily commutative); the category of sheaves of abelian groups on a site. 

\medskip
\noindent
{\bf Definition.} Let $F: \mathsf A\to \mathsf B$ be a functor between two abelian categories. We say that $F$ is {\em exact} (resp. {\em left-exact}) if for every exact sequence $$0\to X'\to X\to X''\to 0$$ in $\mathsf A$, the induced sequence $$0\to F(X')\to F(X)\to F(X'')\to 0\;\;\;\; {\rm (resp.}\;\, 0\to F(X')\to F(X)\to F(X'')\;{\rm )}$$ is exact.

\medskip
\noindent
{\em Comment.} A left-exact functor is exact if and only if it carries surjections to surjections. 

\medskip
\noindent
{\bf Remark.} There is of course also a notion of a {\em right-exact} functor, but we will not use it here.

\medskip
\noindent
{\bf Examples.} 

\medskip
1) Let $\mathsf C$ be a site and let $X$ be an object of $\mathsf C$. Let $\mathsf A$ be the category of sheaves of abelian groups on $\mathsf C$, and let $\mathsf B$ be the category of abelian groups. Let $\Gamma$ be the functor from $\mathsf A$ to $\mathsf B$ that sends a sheaf $\sch F$ to the group $\sch F(X)$. 

By the very definition of an exact sequence of sheaves, $\Gamma$ is left-exact. It is not exact in general: take for $\mathsf C$ the site of open subsets of $\CC$, for $\sch O$ the sheaf of holomorphic functions and for $X$ the open subset $\CC^*$ of $\CC$. Let $d: \sch O\to \sch O$ be the derivation map.  As we have already mentioned, $d$ is surjective but the induced map $\sch O(\CC^*)\to \sch O(\CC^*)$ is not onto ($1/z$ doesn't belong to its image); hence $\Gamma$ is not exact.

\medskip
2) Let $X$ be a Hausdorff topological space, let $\mathsf A$ be the category of sheaves of abelian groups on $X$, and let $\mathsf B$ be the category of abelian groups. Let $\Gamma_c$ be the functor from $\mathsf A$ to $\mathsf B$ that sends a sheaf $\sch F$ to the subgroup $\sch F(X)$ that consists in sections with compact support. By the very definition of an exact sequence of sheaves, $\Gamma_c$ is left-exact, but not exact in general. 

As an example, let $X$ be the open interval $]0;1[$, let $\sch F$ be the sheaf of real analytic functions on $X$ and let $\sch G$ be the sheaf on $X$ that sends an open subset $U$ to $\RR$ if $\frac 1 2\in U$ and to $0$ otherwise; evaluation at $\frac 1 2$ induces a natural map $\sch F\to \sch G$ which is easily seen to be surjective. Nevertheless, $\Gamma_c(\sch F)\to \Gamma_c(\sch G)=0\to \RR$ is not onto; hence $\Gamma_c$ is not exact. 

\medskip
3) Let $G$ be a group. Let $\mathsf A$ be the category of $G$-modules, that is, of abelian groups endowed with an action of $G$ (or, equivalently, of modules over the ring $\ZZ[G]$), and let $\mathsf B$ be the category of abelian groups. Let $F$ be the functor from $\mathsf A$ to $\mathsf B$ that send a $G$-module $M$ to its subgroup $M^G$ consisting of $G$-invariant elements. It is an easy exercise to see that $F$ is left-exact. Let us show by a counter-example that it is not exact in general. 

\medskip
Let us take for $G$ the group $\ZZ/2\ZZ$, acting on $\CC^*$ through complex conjugation. The morphism $\CC^*\stackrel{z\mapsto z^2}\longrightarrow \CC^*$ is equivariant and surjective; but after taking the $G$-invariants one gets $z\mapsto z^2$ from $\RR^*$ to $\RR^*$, which is not onto anymore. 

\medskip
4) There is a topological avatar of the preceeding example, which will play a central role in \'etale cohomology: let $G$ be a {\em profinite group} (e.g. $G={\rm Gal}\;\overline k /k$, for $k$ a field and $\overline k $ an algebraic closure of $k$). Let $\mathsf A$ be the category of {\em discrete} $G$-modules, that is, of $G$-modules such that the stabilizer of every element is an {\em open} subgroup of $G$; it is an abelian category. Let $\mathsf B$ be the category of abelian groups, and let $F$ be the functor from $\mathsf A$ to $\mathsf B$ that sends a $G$-module $M$ to its subgroup $M^G$ consisting of $G$-invariant elements. It is left-exact, but not exact in general (the above counter-example fits in that setting: as $\ZZ/2\ZZ$ is finite, it is profinite and every $\ZZ/2\ZZ$-module is discrete).

\medskip
\noindent
{\bf The notion of a derived functor.} Derived functors of a given left-exact functor are a way to encode its lack of exactness; let us be now more precise. 

\medskip
Let $F:\mathsf A\to \mathsf B$ be a left-exact functor between abelian categories. A collection $(\R^iF)_{i\geq 0}$ of functors from $\mathsf A$ to $\mathsf B$ is said to be a {\em collection of derived functors of $F$} if the following hold:

\medskip
$\bullet$ $\R^0F=F$ ; 

$\bullet$ to every short exact sequence $0\to X'\to X\to X''\to 0$ of objects of $\mathsf A$ one can associate in a functorial way a long exact sequence $$0\to F(X')\to F(X)\to F(X'')\to \R^1F(X')\to \R^1F(X)\to \R^1F(X'')\to\R^2F(X')\to\ldots$$ $$\ldots\to \R^iF(X')\to \R^iF(X)\to \R^iF(X'')\to \R^{i+1}F(X')\to\ldots$$

$\bullet$ the collection $(\R^iF)$ is universal for the preceeding properties, in some straightforward sense we won't explicit. 

\medskip
If such a collection exists, it is unique up to a canonical system of isomorphisms of functors; hence if it is the case, we will speak of {\em the} collection $(\R^iF)$ of the derived functors of $F$, and we will say that $\R^iF$ is {\em the} $i$-th derived functor of $F$. 

\medskip
\noindent
{\bf A stupid example.} If $F$  is an exact functor, then it does admit derived functors: one has $\R^0F=F$ and $\R^iF=0$ for positive $i$. 

\medskip
\noindent
{\bf Conditions for the existence of derived functors.} Under a rather technical condition on the abelian category $\mathsf A$, namely, if it has {\em enough injective objects}, {\bf \em every left-exact functor from $\mathsf A$ to another abelian category $\mathsf B$ does admit derived functors.}

\medskip
Let us now make three remarks. 

\medskip
1) The aforementioned condition is always fulfilled by the following categories we have already encountered: the category of sheaves of abelian groups on a site; the category of modules over a ring, including that of $G$-modules for $G$ a group; the category of {\em discrete} $G$-modules, for $G$ a profinite group\footnote{The latter category is in fact (equivalent to) the category of sheaves of abelian groups on a suitable site.}. Hence on those categories, derived functors will always exist. 

2) The proof of the existence of derived functors (under the above condition) proceeds by {\em building} them, but in a highly abstract way which is usually not helpful at all when one needs to actually compute them. 

3) If $G$ is a (profinite) group, the derived functors of $M\mapsto M^G$, with $M$ going through the category of (discrete) $G$-modules, can be computed explicitely using (continuous) $G$-cochains, cocycles and coboundaries with coefficients in $M$. They will be denoted by $\H^\bullet(G,.)$. 

\medskip
{\small 
\noindent
{\bf About the conditions of having 'enough injective objects'.}  Let $\mathsf A$ be an abelian category. An object $X$ of $\mathsf A$ is said to be {\em injective} if every injective arrow from $X$ to an object of $\mathsf A$ admits a retraction or, otherwise said, identifies $X$ with a direct summand of the target.  As an example, note that in the category of vector spaces over a given field, every object is injective. 

That is not the case for that of abelian groups;  indeed, $\ZZ$ is not an injective object of this category: simply note that the injection $n\mapsto 2n$ from $\ZZ$ into himself doesn't admit a retraction -- indeed, $2\ZZ$ is not a direct summand of $\ZZ$. 

In fact, one can prove (exercise !) that an abelian group is injective if and only if it is divisible. 

\medskip
An abelian category $\mathsf A$ is said to have {\em enough injectives} if every object of $\mathsf A$ admits an injective arrow into an injective object. The category of abelian groups has enough injectives: it is not difficult (again, exercise !) to embed a given abelian group in a divisible one ($\ZZ\hookrightarrow \QQ$ is an example of such an embedding). This fact plays a crucial role in the proof that the category of modules over a given ring and the category of sheaves of abelian groups on a site have enough injectives. 
}

\subsection*{The cohomology groups of a sheaf}

Let us now fix a site $\mathsf C$. Let $X$ be an object of $\mathsf C$; let $\mathsf A$ be the category of sheaves of abelian groups on $\mathsf C$ and let $\mathsf B$ be the category of abelian groups. As we have seen, the functor from $\mathsf A$ to $\mathsf B$ that sends $\sch F$ to $\sch F(X)$ is left-exact; it thus admits derived functors (see above). 

\medskip
\noindent
{\bf Notation and definition.} The functor $\sch F\mapsto \sch F(X)$ will also be denoted $\sch F\mapsto \H^0(X,\sch F)$. For very $i>0$, its $i$'th derived functor will be denoted by $\sch F\mapsto \H^i(X,\sch F)$. We will say that $\H^i(X,\sch F)$ is {\em the $i$-th cohomology group of $X$ with coefficients $\sch F$}, or {\em the $i$-th cohomology group of $\sch F$ on $X$.}

\medskip
\noindent
{\bf Links with other cohomology theories.} Let us assume that $\mathsf C$ is the site of open subsets of a given topological space $X$. For every abelian group $A$, remind that $\underline A$ denotes the sheaf of locally constant $A$-valued functions on $X$. 

\medskip
Then if $X$ is a 'reasonable' topological space (e.g. a metrizable locally finite CW-complex) it can be shown that for every $i$ the group $\H^i(X,\underline A)$ in the above sense does coincide with the $i$-th group of {\em singular} cohomology of $X$ with coefficients $A$. 

\medskip

\medskip
If $Y$ is a closed subset of $X$ (resp. if $X$ is Hausdorff), one can also derive the functor of global sections with support in $Y$ (resp. compact support). The resulting groups are denoted by $\H^\bullet_Y(X,.)$ (resp. $\H^\bullet_c(X,.)$) and called the {\em cohomology groups with support in $Y$ (resp. with compact supports)} of the sheaf involved. Again, under reasonable assumptions on $X$ and $Y$ (resp. on $X$) the group $\H^i _Y(X,\underline A)$ (resp. $\H_c^i (X,\underline A)$) in that sense  coincide with the $i$-th group of {\em singular} cohomology of $X$ with support in $Y$ (resp. with compact support) and with coefficients $\mathsf A$.

\medskip
\noindent
{\bf Some words about the restriction.} Let us still assume that $\mathsf C$ is the site of open subsets of a topological space $X$, let $U$ be an open subset of $X$ and let $V$ be an open subset of $U$; let $\sch F$ be a sheaf of abelian groups on $\mathsf C$. It defines by restriction a sheaf $\sch F_{|U}$ on the topological space $U$. 

\medskip
Now one can define for every $i$, in the above sense, the groups $\H^i(V,\sch F)$ and $\H^i(V,\sch F_{|U})$: for the first one, one computes the $i$-th derived functor of $\sch G\mapsto \sch G(V)$ for $\sch G$ going through the category of abelian sheaves {\em on $X$}, and one applies it to $\sch F$; for the second one, one computes the $i$-th derived functor of $\sch G\mapsto \sch G(V)$ for $\sch G$ going through the category of abelian sheaves {\em on $U$}, and one applies it to $\sch F_{|U}$. 

\medskip
But it turns out that there is no reason to worry about it: both groups $\H^i(V,\sch F)$ and $\H^i(V,\sch F_{|U})$ are canonically isomorphic; hence we will most of the time only use the notation $\H^i(V,\sch F)$. 

\medskip
{\small
\noindent
{\bf About the restriction in a more general setting.} The phenomenon we have just mentioned extends to every site $\mathsf C$, as follows. 

If $X$ is an object of $\mathsf C$, let $\mathsf C/X$ be the so-called category {\em of objects of $\mathsf C$ over $X$}, which is defined in the following way: 

\medskip
$\bullet$ its objects are couples $(Y,f)$ where $Y\in \mathsf C$ and where $f$ is an arrow from $Y$ to $X$;

$\bullet$ if $(Y,f)$ and $(Z,g)$ are two objects of $\mathsf C/X$, an arrow from $(Y,f)$ to $(Z,g)$ is nothing but a morphism $h:Y\to Z$ such that $g\circ h=f$. 

\medskip
The category $\mathsf C/X$ inherits a Grothendieck topology: if $(Y,f)$ is an object of $\mathsf C/X$, a family $((Y_i,f_i)\to (Y,f))$ of arrows in $\mathsf C/X$ is said to be a covering if and only if $(Y_i\to Y)$ is a covering. If $\sch F$ is a sheaf on $\mathsf C$, the assignment $(Y,f)\mapsto \sch F(Y)$ defines a sheaf on $\mathsf C/X$, which is denoted by $\sch F_{|X}$. 

\medskip
Now if $(Y,f)\in \mathsf C/X$, one can define for every integer $i$ the groups  $\H^i(Y,\sch F)$ and $\H^i((Y,f),\sch F_{|X})$: for the first one, one computes the $i$-th derived functor of $\sch G\mapsto \sch G(Y)$ for $\sch G$ going through the category of abelian sheaves {\em on $\mathsf C$}, and one applies it to $\sch F$; for the second one, one computes the $i$-th derived functor of $\sch G\mapsto \sch G(Y,f)$ for $\sch G$ going through the category of abelian sheaves {\em on $\mathsf C/X$}, and one applies it to $\sch F_{|X}$. 

\medskip
But it turns out that there is no reason to worry about it: both groups $\H^i(Y,\sch F)$ and $\H^i((Y,f),\sch F_{|X})$ are canonically isomorphic; hence we will most of the time only use the notation $\H^i(Y,\sch F)$. }

\medskip
\noindent
{\bf Cup-products}. If $\sch F$ and $\sch G$ are two sheaves of abelian groups on a site $\mathsf C$, one has for every $(i,j)$ and every $X\in \mathsf C$ a natural {\em cup-product pairing} $$\cup : \H^i(X,\sch F)\otimes\H^j(X,\sch G)\to \H^{i+j}(X,\sch F\otimes_{\underline \ZZ}\sch G).$$ It is 'graded-commutative' : if $h\in \H^i(X,\sch F)$ and $h'\in \H^j(X,\sch G)$ then one has $h'\cup h=(-1)^{ij}h\cup h'$. 

\subsection*{What is \'etale cohomology ?}

In both scheme-theoretic and Berkovich's settings, one will define the notion of an {\em \'etale map.} 

Now, for $X$ being a scheme (resp. a Berkovich analytic space) we will consider the site $X\etal$ defined as follows: 

\medskip
$\bullet$ its objects are couples  $(U,f)$ where $U$ is a scheme (resp. an analytic space) and $f: U\to X$ an \'etale map (though we often will simply denote such an object by $U$, one has to keep in mind that $f$ is part of the data);

$\bullet$ a morphism from $(V,g)$ to $(U,f)$ is a map $h: U\to V$ such that $f\circ h=g$; {\em such a map is automatically itself \'etale.}

$\bullet$ a family $(g_i : U_i\to U)$ is a covering if and only if it is {\em set-theoretically} a covering, that is, if and only if $U=\bigcup g_i(U_i)$. 

\medskip
We will call $X\etal$ the {\em \'etale site of $X$}, and its topology will be the {\em \'etale topology} on $X$; an {\em \'etale sheaf  on $X$} will be a sheaf on $X\etal$. 

The {\em \'etale cohomology} will then simply denote the cohomology theory of  \'etale sheaves.

\section{\'Etale morphisms of schemes} \label{etmorsch}

General reference: \cite{sga1} or \cite{mil}. 

\subsection*{Motivation, definition and first properties}

Let us begin by some vocabulary. If $Y$ and $X$ are two topological spaces (resp. two complex analytic spaces) and if $y\in Y$ we will say that a continuous (resp. holomorphic) map $f:Y\to X$ is {\em a local homemorphism at $y$} (resp. {\em a local isomorphism at $y$}) if there exists an open neighborhood $V$ of $y$ and an open neighborhood $U$ of $f(y)$ such that $f$ induces an homeomorphism (resp. a holomorphic diffeomorphism) $V\simeq U$. 

\medskip
\noindent
{\bf What is an \'etale map ? } The rough idea one should keep in mind is the following: {\em the fact for a morphism of schemes to be \'etale at a point of the source is an algebraic analog of the fact for a continuous map between topological spaces to be a local homeomorphism at a point of the source.} 

\medskip
Before giving a precise definition, let us give an example of a local homeomorphism coming from differential geometry which we will essentially mimic in the scheme-theoretic context. 
Let $X$ be a $C^{\infty}$-manifold, let $n$ be an integer and let $P_1,\ldots, P_n$ be polynomials in $n$ variables $T_1,\ldots,T_n$ with coefficients in $C^\infty(X)$; let us denote by $Y$ the subset of $X\times\RR^n$ which consists in points $(x,t_1,\ldots,t_n)$ such that $P_i(x,t_1,\ldots,t_n)=0$ for every $i$, and by $f: Y\to X$ the map induced by the first projection. 

Let $y$ be a point of $Y$ at which the Jacobian determinant $\left|\displaystyle{\frac{\partial P_i}{\partial T_j}}\right|_{i,j}$ doesn't vanish; it follows from the implicit function theorem that $f$ is a local homeomorphism at $y$. 

\medskip
\noindent
{\bf Definition.} Let $\phi:Y\to X$ be a morphism of schemes, let $y$ be a point of $Y$ and let $x$ be its image. We will say that $\phi$ is {\em \'etale} at $y$ if there exists an affine neighborhood $U=\spec A$ of $x$ in $X$, a neighborhood $V$ of $y$ in $f^{-1}(U)$, polynomials $P_1,\ldots,P_n$ in $A[T_1,\ldots,T_n]$, and a commutative diagram $$\diagram & \spec A[T_1,\ldots,T_n]/(P_1,\ldots,P_n)\dto \\ V\urto^{\psi} \rto & \spec A\enddiagram$$ with $\psi$ an open immersion such that $\left|\displaystyle{\frac{\partial P_i}{\partial T_j}}\right|_{i,j}$ doesn't vanish on $\psi(V)$. 

\medskip
We will say that $\phi$ is {\em \'etale} if it is \'etale at every point of $Y$. 

\medskip
\noindent
{\bf Some comments.}

\medskip

$\bullet$ It follows from the definition that the \'etale locus of a given map is an open subset of its source. 

$\bullet$ If $\phi : Y\to X$ is a morphism which is \'etale at some point $y$ of $Y$, one can always find a local presentation of $\phi$ as above {\em with $n=1$} (and then one has only to ensure that the partial derivative $\partial P/\partial T$ doesn't vanish on $\psi(V)$). 

$\bullet$ There is another definition of \'etaleness, which is more technically involved: a map $Y\to X$ is \'etale at $y$ if and only if it is {\em flat} and {\em unramified} at $y$.  

$\bullet$ We are now going to give some basic properties of \'etale maps; for the sake of simplicity, we restrict ourselves to {\em global} \'etale morphisms; but most of those results also have a local counterpart. 

\medskip
\noindent
{\bf First properties of \'etale maps}. Every \'etale map is open. If $Z\to Y$ and $Y\to X$ are \'etale, then the composite map $Z\to X$ is \'etale too. If $Y\to X$ is \'etale and if $X'\to X$ is {\em any} morphism, then $Y\times_XX'\to X'$ is \'etale. If $Y\to X$ and $Z\to X$ are \'etale, then {\em any} morphism $Y\to Z$ making the diagram $$\diagram Y \rto \drto &Z\dto \\& X\enddiagram$$ commute is \'etale. 

If you find that the latter property sounds a bit strange, think to the analogous fact for local homeomorphisms in topology -- which is obvious. 

\medskip
\noindent
{\bf Examples of \'etale morphisms.}

\medskip
\begin{itemize}

\item[1)] Every open immersion is \'etale. 

\item[2)] If $Y\to X$ is a morphism between two complex algebraic varieties ({\em i.e} between two schemes of finite type over $\CC$) then it is \'etale if and only if the induced holomorphic map $Y(\CC)\to X(\CC)$ is a local isomorphism on $Y(\CC)$. 

\item[3)] We are now going to give three fundamental examples of {\em finite} \'etale maps (such a map is also called an \'etale cover -- this is the analog of the notion of a topological cover, see further the discussion about the \'etale fundamental group). 

\medskip
\begin{itemize}
\item[$\bullet$] Let $A$ be a commutative ring, let $f\in A^*$, and let $n$ be an integer which is invertible in $A$. The finite morphism $$\spec A[T]/(T^n-f)\to \spec A$$ is \'etale; this follows from the fact that $\displaystyle{\frac{\partial(T^n-f)}{\partial T}}$ doesn't vanish on $\spec A[T]/(T^n-f)$. Indeed, for every point $x$ of the latter one has $T^n(x)=f(x)\neq 0$, whence $nT^{n-1}(x)\neq 0$ (we have used the fact that both $f$ and $n$ are invertible in $A$). 

\medskip
Note that in the particular case where $A=k[X,X^{-1}]$ for $k$ a field in which $n\neq 0$, and where $f=X$, the above finite \'etale map is nothing but $\diagram {\Bbb G}_m\rto ^{x\mapsto x^n}&{\Bbb G}_m\enddiagram$. 

\medskip
\item[$\bullet$] Let $p$ be a prime number, let $A$ be an $\FF_p$-algebra and let $f\in A$. The  finite morphism $$\spec A[T]/(T^p-T-f)\to \spec A$$ is \'etale; this follows from the fact that $\displaystyle{\frac{\partial(T^p-T-f)}{\partial T}}=1$, hence is invertible on $\spec A[T]/(T^p-T-f)$. 

\medskip
Note that in the particular case where $A=k[X]$ for $k$ a field of char. $p$, and where $f=X$, the above finite \'etale map is nothing but $\Aff^1_k\stackrel{x\mapsto x^p-x} {\longrightarrow} \Aff^1_k$. Hence it appears as a {\em non-trivial \'etale cover of the affine line} (of degree $p$, by itself). Note that this example is 'geometric' and not 'arithmetic' in the sense that it holds for {\em every} field of char. $p$, even algebraically closed. 

\medskip
We thus see that {\em the affine line over an algebraically closed field of char. $p$ is never simply connected}. 

\medskip
\item[$\bullet$] Let $k$ be a field. A finite extension $L$ of $k$ is said to be {\em separable} if it is isomorphic to $k[X]/P$ for some irreducible $P\in k[X]$ with $P'\neq 0$ (this is equivalent to the fact that $(P,P')=1$). If char. $k=0$, if $k$ is algebraically closed or if $k$ is finite, every finite extension of $k$ is separable. A finite $k$-algebra is said to be separable if it is a product of finitely many finite separable extensions of $k$. 

One can prove that if $X$ is a scheme of finite type over $k$, it is \'etale over $k$ if and only if it is isomorphic to the spectrum of a  finite separable $k$-algebra (and $X$ is then finite \'etale over $k$). 

Note that if $k$ is algebraically closed, a finite separable $k$-algebra is an algebra isomorphic to $k^n$ for some $n$ ; hence a $k$-scheme of finite type is \'etale over $k$ if and only if it is a disjoint sum of finitely many copies of $\spec k$. 

The latter assertion is a very simple illustration of the following principle: from the \'etale viewpoint, the spectrum of an algebraically (and, more generally, separably) closed field is analog to a one-point space in topology; this is definitively false for the spectrum of an arbitrary field: for example, as $\CC$ is a separable extension of $\RR$, the map $\spec \CC\to \spec \RR$ is a {\em non-trivial \'etale cover of degree 2} of $\spec \RR$.

\end{itemize}

\end{itemize}

\medskip
\noindent
{\bf The fibers of an \'etale map.} Let $\phi : Y\to X$ be an \'etale morphism of schemes, and let $x\in X$. The fiber $\phi^{-1}(x)$ is discrete, and finite soon as $\phi$ is of finite type ({\em e.g.} $Y$ and $X$ are two algebraic varieties over a given field); scheme-theoretically one can write $\phi^{-1}(x)=\coprod \spec k_i$, with each $k_i$ finite separable over $\kappa(x)$. 

\medskip
\noindent
{\bf The \'etale topology on a scheme}. 

Let $X$ be a scheme. We have defined what an \'etale morphism of schemes is. This allows us to define the \'etale site $X\etal$ and to speak about the \'etale topology on $X$ (see just before section \ref{etmorsch}); we now want to give an example of a property which is satisfied locally for the \'etale topology, but not in general for the Zariski one. 

\medskip
Let $Y\to X$ be a morphism of schemes, and $y\in Y$. We will say that $Y\to X$ is {\em smooth at $y$} if there exists an integer $n$, an open neighborhood $U$ of $y$ and a factorization of $U\to X$ through an \'etale map $U\to \Aff^n_X$; note that if $Y\to X$ is \'etale at $y$, it is smooth at $y$ (take $n=0$ and $U$ beeing the \'etale locus of $Y\to X$). 

\medskip
Now let $k$ be a field and let $X$ be a non-empty, smooth $k$-scheme. By definition of smoothness, there exists a non-empty open subset $U$ of $X$
 and \'etale map $U\to \Aff^n_k$ for some $n$. As $U\to \Aff^n_k$ is \'etale, its image is a non-empty open subset $V$ of $\Aff^n_k$. 
 
 \medskip
 Let $k^s$ be a separable closure of $k$; it is an infinite field, therefore $V$ admits a $k^s$-point; but this exactly means that there exists $x\in V$ such that $\kappa(x)$ embeds into $k^s$, hence is finite separable over $k$. 
 
 By the very definition of $V$, there exists $y\in U$ whose image on $\Aff^n_k$ is $x$. As $V\to U$ is \'etale, $\kappa(y)$ is finite separable over $\kappa(x)$, hence over $k$.

 \medskip
 We therefore have shown that $X$ has a point whose residue field is separable over $k$ ; let us rephrase this property in two different ways which are tautologically equivalent (it follows from the universal property of the fiber product). 
 
 \medskip
 {\em First way.} There exists a finite separable extension $F$ of $k$ and a commutative diagram $$\diagram & X\dto \\ \spec F \urto \rto  & \spec k\enddiagram$$
 
 {\em Second way}. There exists a finite separable extension $F$ of $k$ and a section of the morphism $(X\times_{\spec k}\spec F)\to \spec F$. 
 
 \medskip
 For very finite separable extension $F$ of $k$, the map $\spec F\to \spec k$ is \'etale; as it is trivially surjective ($\spec k$ is a point !), it is an \'etale covering. We therefore have shown that every non-empty smooth scheme over a field $k$ has a section {\em locally for the \'etale topology on $\spec k$.}
 
 \medskip
 This fact admits the following generalization, whose proof is a little bit more technical, but essentially no more difficult, than the preceding one: {\em if $\phi: Y\to X$ is a smooth, {\em surjective} morphism between two schemes, then it does admit a section locally for the \'etale topology on $X$.}
 
 \medskip
 \noindent
 {\bf Some comments.} 
 
 \medskip
 1) The above claim is definitely false if one restricts to Zariski topology, which is too coarse. Let us give two counter-examples, the first one being 'arithmetic' ({\em i.e.} related to the fact that the ground field is not algebraically closed) and the second one being 'geometric' ({\em i.e.} over an algebraically closed field). 
 
 \medskip
 {\em The arithmetic example.} Let $X$ be the non-trivial real conic, that is, the projective $\RR$-scheme defined in homogeneous coordinates by the equation $$T_0^2+T_1^2+T_2^2=0 ; $$it can be easily shown that $X$ is smooth, but $X\to \spec \RR$ doesn't admit any section, because $X(\RR)=\emptyset$ (as $\spec \RR$ is a point, to have a section locally for the Zariski topology on $\spec \RR$ just means to have a section). Of course, $X(\CC)\neq \emptyset$. 
 
 \medskip
 
 {\em The geometric example.} Set $Y=\spec \CC[U,V,S,T]/(SU^2+TV^2-1)$, let $\phi: Y\to \Aff^2_{\CC}$ be the morphism defined by $S$ and $T$, and let 
 $X$ be the complement of the origin in $\Aff^2_{\CC}$. One can check that $\phi$ is smooth over $X$ ; but it doesn't admit sections locally on $X$. 
 
 Indeed, one can even prove that if $X'$ is any non-empty open subset of $X$, then $\phi$ doesn't admit a section over $X'$: if it did, the restriction of this section to the generic fiber $Y_\eta$ of $\phi$ would define a $\CC(S,T)$-point on $Y_\eta$, that is, a solution of the equation $SU^2+TV^2=1$ in $\CC(S,T)$ (variables are $U$ and $V$); but (exercise) there are no such solutions. 
 
 \medskip
 Let us now show concretely that $\phi$ has locally a section on $X$ for the \'etale topology. For that purpose, let us denote by $\psi_1$ the map $$\Aff^2_{\CC}=\spec \CC[S',T]\to \Aff^2_{\CC}=\spec \CC[S,T]$$ that is induced by $S\mapsto  (S')^2$ and $T\to T$ (at the level of points, it is simply $(s',t)\mapsto ((s')^2, t)$), and by $\psi_2$ the map $$\Aff^2_{\CC}=\spec \CC[S,T']\to \Aff^2_{\CC}=\spec \CC[S,T]$$ that is induced by $S\mapsto  S$ and $T\to (T')^2$ (at the level of points, it is simply $(s,t')\mapsto (s, (t')^2)$).

 Let $X_1$ (resp. $X_2$) be the open affine subscheme of $\Aff^2_{\CC}$ defined as the invertibility locus of $S$ (resp. $T$). It is easily seen that $\psi_i^{-1}(X_i)\to X_i$ is a finite \'etale cover of degree 2 for $i=1,2$; as $X=X_1\cup X_2$, the family $(\psi_1^{-1}(X_1)\to X,\psi_2^{-1}(X_2)\to X)$ is an \'etale covering of $X$. 
 
 \medskip
 Now for every $i$ the map $Y\times_{\Aff^2_{\CC}} \psi^{-1}(X_i)\to \psi^{-1}(X_i)$ admits a section: 
 
 \hspace{.5cm}$\bullet$ if $i=1$ one can consider the map defined as the level of rings by the formulas $U\mapsto (1/S'), V\mapsto 0$ (at the level of points, this map corresponds to $(s',t)\mapsto (1/s', 0, s',t)$ ); 
 
 \hspace{.5cm} $\bullet$ if $i=2$ one can consider the map defined as the level of rings by the formulas $U\mapsto 0, V\mapsto 1/T'$ (at the level of points, this map corresponds to $(s,t')\mapsto (0, 1/t', s,t')$ ). 
 
 \medskip
 2) In the complex analytic case, it is obvious that a smooth surjective map admits sections locally on the target (because such a map is locally on the source and on the target isomorphic to a projection $X\times \DD^n\to X$ for some $n\geq 0$, where $\DD$ is the open unit disc); hence we see that the \'etale topology is, for some purposes, a good algebraic analog of the classical topology over $\CC$ -- at least far better than the Zariski topology. 
 
\subsection*{The fundamental group of a scheme}

As we will see, Grothendieck's theory of the fundamental group provides, through the notion of a finite \'etale map, a very natural frame containing both Galois theory and the theory of (finite) {\em topological} coverings of complex algebraic varieties, which may help to get a geometric intuition of some {\em a priori} purely algebraic field-theoretic notions. 

\medskip
In consistency with the fact that the \'etale avatar of a one-point space is supposed to be the spectrum of a separably closed field, one defines a {\em geometric point} on a scheme $X$ as a morphism $\spec F\to X$, for $F$ a separably closed field; equivalently, such a point can be defined as a datum of a (true) point $x$ of the scheme $X$, together with an embedding $\kappa(x)\hookrightarrow F$ for some separably closed field $F$. 
 
\medskip
Let $X$ and $Y$ be two connected, non-empty schemes and let $Y\to X$ be finite and \'etale. It can be shown that $G:=\mathsf{Aut}\;Y/X$ is a finite group. 
 
 One says that $Y\to X$ is {\em Galois} if $X$ is {\em scheme-theoretically} equal to the quotient $Y/G$; it means exactly that every arrow $Y\to Z$ which is G-invariant goes in a unique way through an arrow $X\to Z$, and if $Y=\spec A$, this is equivalent to the fact that $X=\spec A^G$. 
 
 As an example, if $k$ is a field and if $L$ is a finite separable extension, then the finite \'etale map $\spec L\to \spec k$ is Galois if and only if $L$ is a Galois extension of $k$ in the usual sense.

 \medskip
 \medskip
 Now, let us choose a {\em geometric point} $\overline x$ on $X$ (one should think to it as a 'basepoint'). One can associate to it a profinite group $\pi_1(X,\overline x)$ which is called the {\em fundamental group} of $(X,\overline x)$.  
 
 \medskip
 The group $\pi_1(X,\overline x)$ encodes the theory of finite \'etale covers of $X$. This quite vague claim can be given a precise meaning, involving an equivalence of categories, which we won't explain in full detail here; let us nevertheless give some precisions. 
  
 \medskip
 $\bullet $ To every open subgroup $U$ of $\pi_1(X,\overline x)$ is associated a connected, non-empty finite \'etale cover $\mathsf {Cov}\;U$ of $X$.
 
 $\bullet$ If $U$ and $V$ are two open subgroups of  $\pi_1(X,\overline x)$ with $U\subset V$, the arrow $\mathsf {Cov}\;U \to X$ goes through a natural finite \'etale cover $\mathsf {Cov}\;U\to \mathsf {Cov}\;V$. 
 
 $\bullet$ If $U$ is an open subgroup $U$ of $\pi_1(X,\overline x)$ then $\mathsf {Cov}\;U$ is Galois if and only if $U\lhd \pi_1(X,\overline x)$, and if it is the case there is a natural isomorphism $\mathsf{Aut}\;(\mathsf{Cov}\;U/X)\simeq \pi_1(X,\overline x)/U$.
 
 $\bullet$ Every connected, non-empty, finite \'etale cover is isomorphic to $\mathsf{Cov}\;U$ for some $U$. 
 
 $\bullet$ $\pi_1(X,\overline x)$ is naturally isomorphic to $$\lim_{\stackrel \leftarrow {U\lhd \pi_1(X,\overline x)}}\mathsf{Aut}\;(\mathsf{Cov}\;U/X).$$
 
 \medskip
 Let us mention that if $\overline y$ is another geometric point of $X$, then $\pi_1(X,\overline y)$ is {\em non-canonically} isomorphic to $\pi_1(X,\overline x)$. Therefore it may happen that one speaks of {\em the} fundamental group of $X$, without mentioning any geometric point  (and that one notes it $\pi_1(X)$)-- but this is a little bit abusive and sometimes dangerous, exactly like speaking of the fundamental group of a connected, non-empty topological space, without having pointed it. 
 
 \medskip
 \noindent
 {\bf Example.} Suppose $X=\spec k$, with $k$ a field, and $\overline x=\spec k^s$ for some separable closure $k^s$ of $x$. Then $\pi_1(X,\overline x)$ is nothing but $\mathsf{Gal}\;(k^s/k)$, and what is written above is (part of) Galois correspondence. 
 
 \medskip
 \noindent
 {\bf Remark.} We  see that choosing a separable closure of a given field is, in some sense, analog to choosing a basepoint on a (connected, non-empty) topological space. 
 
 \medskip
 \noindent
 {\bf Functoriality.} Let $Y\to X$ be a morphism between two connected, non-empty schemes, and let $\overline y$ be a geometric point of $Y$. It defines by composition with $Y\to X$ a geometric point $\overline x$ of $X$; the map $Y\to X$ then induces in a natural way a continuous morphism $\pi_1(Y,\overline y)\to \pi_1(X,\overline x)$.

\medskip
\noindent
{\bf Various comparison theorems}. In what follows, we will implicitly use the following fact: let $k$ be an algebraically closed field and let $L$ be an algebraically closed extension of $k$.  If $X$ is a connected, non-empty algebraic variety over $k$, then the $L$-algebraic variety $X_L:=X\times_kL$ is non-empty and connected too; if $k=\CC$, then $X(\CC)$ is also non-empty and connected for the transcendent topology. 

\medskip
Our aim is now to describe the behaviour of the fundamental group when one increases the {\em algebraically closed} ground field, and how the Grothendieck fundamental group of an algebraic variety over $\CC$ and its classical (topological) fundamental group are related. 

\medskip
 \noindent
 {\bf Theorem 1}. {\em Let $k$ be an algebraically closed field and let $X$ be a connected, non-empty algebraic variety over $k$. Let $L$ be an algebraically closed extension of $k$, let $\overline \xi$ be a geometric point of $X_L$ and let $\overline x$ be the geometric point on $X$ defined by composition of $\overline \xi$ with $X_L\to X$. 
 
 Assume that $k$ is of char. $0$ or that the variety $X$ is proper (e.g. projective). The  natural map $\pi_1(X_L,\overline \xi)\to \pi_1(X,\overline x)$ is then an isomorphism. }
 
 \medskip
 \noindent 
 {\bf Some comments}. 
 
 \medskip
 1) In a more concrete way, what theorem 1 says is that if char $k=0$ or if $X$ is proper then {\em every finite \'etale covering of $X_L$ can already be defined over $k$.} 
 
 2) Theorem 1 is definitely false if char. $k=p>0$ and if $X$ is not proper, in which case the fundamental group will become huger and huger by enlarging the algebraically closed ground field; or, in other words, there will be plenty of \'etale covers of $X_L$ not definable over $k$ as soon as $L\neq k$.
 
Let us give an example of such a cover: if $\lambda\in  L-k$, the \'etale cover $$\spec L[X,T]/(T^p-T-\lambda X)\to \spec L[X]=\Aff^1_L$$ {\em can not be defined over $k$.}

3) Let us assume that char. $ k=p>0$; for every profinite group $G$, we denote by $G_{\neq p}$ the largest profinite quotient of $G$ all of whose finite quotients are of prime-to-$p$-order. One can then prove, with the notations of th. 1 and without any properness assumption, that $\pi_1(X_L,\overline \xi)_{\neq p}\to \pi_1(X,\overline x)_{\neq p}$ is an isomorphism. This means that every {\em prime-to-$p$} Galois cover of $X_L$ can be defined over $k$. 

\medskip
Let us come now to the comparison between the \'etale and topological fundamental groups. 
 
\medskip
 \noindent
 {\bf Theorem 2}. {\em Let $X$ be a connected, non-empty complex algebraic variety, let $x\in X(\CC)$ and let us denote by $x$ the corresponding geometric point $\spec \CC\to X$ of the scheme $X$. There is a natural isomorphism between $\pi_1(X, x)$ and the profinite completion of $\pi_{1,\rm top}(X(\CC),x)$.}
 
 \medskip
 \noindent
 {\bf Some comments.} In a more concrete way, what theorem 2 says is that {\em every finite topological cover of $X(\CC)$ is algebraizable}; this is what is classically called 'Riemann's existence theorem'.  
 
 \medskip
 \noindent
 {\bf A consequence of th. 1 and 2.} Let $X$ be a connected, non-empty algebraic variety over $\CC$ and let $\sigma$ be any (possibly non-continuous) automorphism of $\CC$. The corresponding twisted variety $X_\sigma$ coincides with $X$ as a scheme; only its structure map to $\spec \CC$ has been changed (through $\sigma$). 
 
 By applying th. 2 to $X$ and $X_\sigma$, and by using the fact that they coincide as schemes, one sees that {\em the groups $\pi_{1,\rm top}(X(\CC))$ and $\pi_{1,\rm top}(X_\sigma(\CC))$ have isomorphic profinite completions.} 
 
 \medskip
 Note that this assertion about profinite completions simply means that $X(\CC)$ and $X_\sigma(\CC)$ 'have the same finite topological covers'.  
 
 \medskip
 \noindent
 {\bf Remak.} Serre has given an example (\cite{sgf}) of a complex algebraic variety $X$ and an automorphism $\sigma$ of $\CC$ such that $\pi_1(X(\CC))$ and $\pi_1(X_\sigma(\CC))$ are {\em not} isomorphic though they have, by the above, the same profinite completions. Roughly speaking, $X(\CC)$ and $X_\sigma(\CC)$ have the same finite topological covers, but not the same general topological covers. 
 
\medskip
\noindent
{\bf The fundamental group of a projective curve.} Let $k$ be an algebraically closed field and let $X$ be a smooth, irreducible, projective $k$-curve. 

\medskip
$a)$ If $X=\PP^1_k$ then $\pi_1(X)=0$. 

$b)$ If $X$ is of positive genus $g$ and if char. $k=0$, then $\pi_1(X)$ admits, {\em as a profinite group}, a presentation $$<a_1,b_1,\ldots,a_g,b_g>/[a_1,b_1][a_2,b_2]\ldots[a_g,b_g].$$ 

$c)$ If $X$ is of positive genus $g$ and if char. $k=p>0$, then $\pi_1(X)_{\neq p}$ admits, {\em as a profinite group with only prime-to-$p$ finite quotients}, a presentation $$<a_1,b_1,\ldots,a_g,b_g>/[a_1,b_1][a_2,b_2]\ldots[a_g,b_g].$$

\medskip
\noindent
{\small 
{\bf About the proofs of the above results.}  Claim a) follows easily from Hurwitz's formula; claim b) is proven by reducing, through th. 1, first of all to the case where $k$ is contained in $\CC$, and then to the case where $k=\CC$, in which one just has to apply th. 2 (and the transcendental standard description of the fundamental group of a Riemann surface). 

\medskip
Claim c) is far more subtle and was the first success of Grothendieck's approach to algebraic geometry. 

{\em First step.} One exhibits a complete, discrete valuation ring $A$ with residue field $k$ and with quotient field $K$ of char. $0$: this is the 'Witt vectors construction'. 

{\em Second step.} One extends $X$ to a flat, smooth, projective $A$-scheme $\sch X$, as follows: first, one makes a formal extension, which is possible because the obstruction to do it lives in a Zariski $\H^2$ group of $X$, which is $0$ by a dimension argument; and then, one algebraizes this formal extension, which is again possible thanks to dimension 1. 

{\em The core of the proof.} Let $\overline K$ be an algebraic closure of $K$. A comparison theorem which is due to Grothendieck says that $$\pi_1(\sch X\times_A \overline K)_{\neq p}\simeq \pi_1(X)_{\neq p}\; ; $$ it is the scheme-theoretic analog of the following classical fact: if $\phi : Y\to X$ is a smooth, surjective and proper map with connected fibers between complex analytic spaces, then $x\mapsto \pi_{1\rm top}(\phi^{-1}(x))$ is locally constant on $X$.

{\em Conclusion.} The required conclusion now follows by applying assertion b) to the curve $\sch X\times_A \overline K$.} 

\section{\'Etale morphisms of analytic spaces}

General reference: \cite{brk2}. 

\medskip
{\em We fix a complete, non-Archimedean field $k$.  }

\medskip
In order to simplify exposition, we have chosen to give the definition of a finite \'etale map before that of a general \'etale map. But fortunately, one can check that \'etale maps which are moreover finite are nothing but the 'finite \'etale' maps defined in a first time: though it is not introduced in the logical way, our terminology is consistent.

\medskip
\noindent
{\bf Definition.} Let $\phi : Y\to X$ be a morphism between two $k$-analytic spaces. We say that $\phi$ is {\em finite \'etale} if for every affinoid domain $U=\sch M(\sch A)$ of $X$, the pre-image $V$ of $U$ on $Y$ is an affinoid domain of $Y$ equal to $\sch M(\sch B)$ for some finite \'etale $\sch A$-algebra $\sch B$. 

\medskip
\noindent
{\bf Remark.} It is sufficient to check the above property on an affinoid G-covering of $X$. 

\medskip
Let us now come to the definition of general \'etale maps. 

\medskip
\noindent
{\bf Definition.} Let $\phi : Y\to X$ be a morphism between two $k$-analytic spaces, let $y\in Y$ and let $x=\phi(y)$. One says that $\phi$ is {\em \'etale at $y$} if there exists an open neighborhood $U$ of $x$ and an open neighborhood $V$ of $y$ such that $\phi$ induces a finite \'etale map (in the above sense) from $V$ to $U$. 

We say that a morphism $Y\to X$ between two $k$-analytic spaces is {\em \'etale} if it is \'etale at every point of $Y$.

\medskip
\noindent
{\bf Remark.} Be aware that one can not, in general, take $V$ being equal to $\phi^{-1}(U)$. 

\medskip
\noindent
{\small
{\bf Zariski topology is too coarse to make the above definition work in scheme theory.} Indeed, let us assume that $p\neq 2$ ; let $Y$ be the invertibility locus of $T-1$ on ${\Bbb G}_m$ and let $\phi$ be the map $Y\to {\Bbb G}_m$ induced by the finite \'etale arrow ${\Bbb G}_m\stackrel{z\mapsto z^2}\longrightarrow {\Bbb G}_m$. Let $y$ be the closed point of $Y$ at which $T+1=0$; its image $x$ is the closed point of ${\Bbb G}_m $ at which $T=1$, and one has $\kappa(y)=\kappa(x)=k$. 

Then though $\phi$ is \'etale at $y$, it is not possible to find an open neighborhood $V$ of $y$ in $Y$ and an open neighborhood $U$ of $x$ in ${\Bbb G}_m$ such that $\phi$ induces a finite \'etale map $V\to U$. Indeed, if it were possible, then by checking what happens at the generic point of the target, one sees that the degree of $V\to U$ would be equal to $2$, contradicting the fact that $\phi^{-1}(x)=\{y\}$ with $\kappa(y)=\kappa(x)$, hence is of degree {\em one} over $\kappa(x)$. 

}

\medskip
\noindent
{\bf First properties.} The \'etale locus of a given map is an open subset of its source; if $\phi: Y\to X$ is a morphism of $k$-analytic spaces, then $\phi$ is \'etale at a point $y\in Y$ if and only if $\phi$ is finite, flat and unramified at $y$. 

\medskip
We are now going to give some basic properties of \'etale maps; for the sake of simplicity, we restrict ourselves to {\em global} \'etale morphisms; but most of those results also have a local counterpart. 

\medskip
Every \'etale map is open and {\em boundaryless}. If $Z\to Y$ and $Y\to X$ are \'etale, then the composite map $Z\to X$ is \'etale too. If $Y\to X$ is \'etale and if $X'\to X$ is {\em any} morphism, then $Y\times_XX'\to X'$ is \'etale. If $Y\to X$ and $Z\to X$ are \'etale, then {\em any} morphism $Y\to Z$ making the diagram $$\diagram Y \rto \drto &Z\dto \\& X\enddiagram$$ commute is \'etale. 

\medskip
\noindent
{\bf Examples.}

\begin{itemize}
\medskip
\item[$\bullet$] Every open immersion is \'etale. Be careful: the embedding of a non-open analytic domain ({\em e.g.} the closed unit ball in an affine space) is {\em not} \'etale because it has non-empty boundary. 

\item[$\bullet$] If $Y\to X$ is a morphism between algebraic $k$-varieties, then $Y\an \to X\an$ is \'etale if and only if $Y\to X$ is \'etale. 

\item[$\bullet$] A very important class of \'etale maps is that of {\em \'etale covers.} A morphism $\phi : Y\to X$ is said to be an \'etale cover if $X$ admits an open covering $(X_i)$ such that for every $i$ one can write $$\phi^{-1}(X_i)=\coprod_{j\in J_i} Y_{i,j}$$ with each $Y_{i,j}$ finite \'etale over $X_i$ (and with the index set $J$ possibly infinite). Let us now give three examples of \'etale covers.

\medskip
\begin{itemize}

\item[$1)$] {\em An example of finite \'etale cover.} Let $r$ and $R$ be two positive real numbers with $0<r<R$, let $n$ be a positive integer which is non-zero in $k$, and let $Y$ (resp. $X$) be the open subset of $\Aff^{1,an}_k$ defined by the conditions $r<|T|<R$ (resp. $r^n<|T|<R^n$). The map $$Y\to X, z\mapsto z^n$$ is a finite \'etale cover of degree $n$. 

\item[$2)$] {\em A topological cover.} Every topological cover of an analytic space inherits a structure of an analytic space -- and becomes an \'etale cover of the target. Let us give now a fundamental example.

Assume that $k$ is algebraically closed and non-trivially valued, and let $E$ be a $k$-elliptic curve with {\em bad} reduction. By a celebrated theorem of Tate (which was written in the rigid-analytic language and was actually the first motivation for developing such a theory), there exist $q\in k^*$ with $|q|<1$ and an isomorphism $E\an\simeq \gm \an/q^\ZZ$. 

The {\em uniformization} arrow $\gm\an\to E\an$ is a topological cover -- this is even a universal cover of $E\an$ --, and in particular it is an \'etale cover. 

\item[$3)$] {\em The logarithm.} Let $\CC_p$ be the completion of an algebraic closure of $\QQ_p$, let $X$ be the open unit $\CC_p$-disc with center 1, and let $\log$ be the usual logarithm defined on $X$ by the formula $$\log (1+z)=z-\frac {z^2}2+\frac{z^3}3-\ldots$$ Then $\log$ induces an infinite \'etale cover $X\to \Aff^{1,an}_{\CC_p}$; it is not a topological cover (indeed, $\Aff^{1,an}_{\CC_p}$ is topologically simply connected).

\end{itemize}

\item[$\bullet$] {\em Some comments on the above examples}. The two last examples show that, countrary to what happens in algebraic geometry, one can exhibit in analytic geometry natural {\em infinite} \'etale covers between connected (and very basic) spaces. 

\item[$\bullet$] {\em The notion of a Galois cover}. As in the algebraic setting, an \'etale cover between analytic spaces is said to be Galois if the target is the quotient of the source (in a suitable meaning, analog to the one we used in scheme-theory) by the automorphism group of the cover. 

Example 1) above is Galois as soon as $k$ is algebraically closed, and its Galois group is then equal to $\mu_n(k)$ -- the latter beeing itself non-canonically isomorphic to $\ZZ/n\ZZ$. 

Example 2) above is Galois, with Galois group $q^\ZZ$, hence isomorphic to $\ZZ$. 

Example 3) above is Galois, with Galois group $\mu_{p^\infty}(\CC_p):=\bigcup\limits_n \mu_{p^n}(\CC_p)$. 
\end{itemize}

\medskip
\noindent
{\bf The fibers of an \'etale map.} Let $\phi : Y\to X$ be an \'etale morphism between two $k$-analytic spaces and let $x\in X$. Its fiber $\phi^{-1}(x)$ is discrete; this is a disjoint (possibly infinite) union $\coprod \sch M(k_i)$ where every $k_i$ is a finite separable extension of $\hres(x)$. 

\medskip
Moreover, if $y$ is a preimage of $x$, the germ of morphism $(Y,y)\to (X,x)$ {\em only depends on the finite separable field extension $\hres(x)\hookrightarrow \hres(y)$.} 

Note a very useful particular case of the latter assertion: {\em the map $\phi$ is a local isomorphism at $y$ if and only if $\hres(x)=\hres(y)$.}

\medskip
If $k$ is algebraically closed and if $x$ is a $k$-point, then $\hres(y)=\hres(x)$ for every preimage $y$ of $x$, hence $\phi$ is a local isomorphism at every pre-image of $x$. 

\medskip
But this is false in general, that is, if $k$ is not algebraically closed or if $x$ is not a $k$-point; let us give a simple, explicit example.  For every $s>0$ denote by $\eta_s$ the point of $\Aff^{1,an}_k$ defined by the semi-norm $\sum a_iT^i\mapsto \max |a_i|s^i$. Then in example 1) above for every $s\in ]r;R[$ the pre-image of $\eta_{s^n}$ is $\{\eta_{s}\}$, and the field extension $\hres(\eta_{s^n})\hookrightarrow \hres(\eta_s)$ induced by our \'etale map is of degree $n$, hence is non-trivial as soon as $n>0$. 

\medskip
\noindent
{\bf The \'etale topology on an analytic space}. 

Let $X$ be an analytic space. We have defined what an \'etale morphism of analytic spaces is. This allows us to define the \'etale site $X\etal$ and to speak about the \'etale topology on $X$ (see just before section \ref{etmorsch}); we now want to give an example of a property which is satisfied locally for the \'etale topology on $X$,  but not in general for its usual topology. 

\medskip
Let $Y\to X$ be a morphism of analytic spaces, and let $y\in Y$. We will say that $Y\to X$ is {\em smooth at $y$} if there exists an integer $n$, an open neighborhood $U$ of $y$ and a factorization of $U\to X$ through an \'etale map $U\to \Aff^n_X$; note that if $Y\to X$ is \'etale at $y$, it is smooth at $y$ (take $n=0$ and $U$ beeing the \'etale locus of $Y\to X$). 

\medskip
As in scheme theory, the following is true: {\em if $\phi: Y\to X$ is a smooth, {\em surjective} morphism between two analytic spaces, then it does admit a section locally for the \'etale topology on $X$.}
 
 \medskip
 \noindent
 {\bf Comment.} The above claim is definitely false if one restricts to the usual topology on $X$, which is too coarse.

\medskip
\noindent
{\bf Fundamental groups of an analytic space.} Let $X$ be a non-empty, connected $k$-analytic space and let $\overline x$ be a geometric point of $x$, that is, a true point $x$ of $X$ together with an embedding of $\hres(x)$ into a separably closed field. To the datum $(X,\overline x)$ one can associate in a functorial way various topological groups, which are in general neither discrete nor profinite, each of which encoding part of the theory of \'etale covers of $X$. 

\medskip
As an example, the theory of the finite \'etale covers (resp. topological \'etale covers, resp. all \'etale covers) of $X$ is encoded in a profinite (resp. discrete, resp. topological) group which is denoted by $\pi_{1,\rm f}(X,\overline x)$ (resp. $\pi_{1,\rm top}(X,\overline x)$, resp. $\pi_1(X,\overline x)$); Yves Andr\'e has introduced another one, which is denoted by $\pi_{1,\rm temp}(X,\overline x)$ and which classifies the so-called {\em tempered} \'etale covers, that is, those which become topological after a finite \'etale base change; this group has turned out to be of deep arithmetical interest. 

\medskip
\noindent
{\em First remark.}  Both profinite completions of $\pi_1(X,\overline x)$ and $\pi_{1,\rm temp}(X,\overline x)$ are naturally isomorphic to $\pi_{1,\rm f}(X,\overline x)$. 

\medskip
\noindent
{\em Second remark.} Each of those groups only depends, {\em up to a non-canonical isomorphism}, on $X$ and not on the geometric point $\overline x$; therefore it may happen that one doesn't mention the latter. 

\medskip
\noindent
{\bf Examples.} 

\medskip
$\bullet$ If $X$ is an algebraic variety over $k$ and if $\overline x$ is a geometric point of $X\an$, then it induces a geometric point $\overline \xi$ of $X$, and there is a natural continuous homomorphism $\pi_{1,\rm f}(X\an, \overline x)\to \pi_1(X,\overline \xi)$. 

If char. $k=0$ then $\pi_{1,\rm f}(X\an, \overline x)\to\pi_1(X,\overline \xi)$ is an isomorphism; if char. $k=p>0$ then it only induces an isomorphism $\pi_{1,\rm f}(X\an, \overline x)_{\neq p}\to \pi_1(X,\overline \xi)_{\neq p}$. Roughly speaking, if char. $k=0$ (resp. char. $k=p>0$) then every finite (resp. finite and prime-to-$p$) Galois cover of $X\an$ is algebraic. 

$\bullet$ Suppose that $k$ is algebraically closed and that $X$ is an open disc over $k$; let $p$ be the characteristic of the residue field $\red k$. If $p=0$ (resp. if $p>0$) then every finite (resp. finite and prime-to-$p$) Galois cover of $X$ is trivial; hence $\pi_{1,\rm f}(X)=0$ (resp.  $\pi_{1,\rm f}(X)_{\neq p}=0$). 

$\bullet$ Suppose that $k$ is algebraically closed and that $X$ is an open annulus over $k$; let $p$ be the characteristic of the residue field $\red k$. If $p=0$ (resp. if $p>0$) then every finite (resp. finite and prime-to-$p$) Galois cover of $X$ is Kummer, that is, as in example 1) above; hence $\pi_{1,\rm f}(X)\simeq \widehat{\ZZ}$ (resp.  $\pi_{1,\rm f}(X)_{\neq p}\simeq \widehat{\ZZ}_{\neq p}=\prod\limits_{\ell \neq p}\ZZ_\ell$). 

$\bullet$ Suppose that $k$ is algebraically closed and let $E$ be an elliptic curve over $k$. If $E$ has good reduction then $$\pi_{1,\rm f}(E\an)\simeq \widehat{\ZZ}^2\;\;; \pi_{1,\rm top}(E\an)=0\;\;\;{\rm and}\; \pi_{1,\rm temp}(E\an)\simeq \widehat{\ZZ}^2\;; $$ if $E$ has bad reduction then $$\pi_{1,\rm f}(E\an)\simeq \widehat{\ZZ}^2\;\;; \pi_{1,\rm top}(E\an)\simeq \ZZ \;\;\;{\rm and}\; \pi_{1,\rm temp}(E\an)\simeq \widehat{\ZZ}\times \ZZ\;. $$

$\bullet$ The topological group $\pi_1(\PP^{1,an}_{\CC_p})$ is huge. 

\section{\'Etale sheaves and \'etale cohomology on schemes}

General reference: \cite{sga42}, \cite{sga43},  \cite{sga45}, \cite{sga5} (for the $\ell$-adic formalism) and \cite{mil}. 

\medskip
Let us begin with a remark: as open immersions are \'etale, any open covering of a scheme $X$ is an \'etale covering of $X$; therefore every sheaf on $X\etal$ is, when restricted to the category of open subsets of $X$, a usual sheaf. 

\subsection*{Examples of sheaves on $X\etal$}

\medskip
\noindent
{\bf The constant sheaves.} Let $A$ be a set. The constant sheaf $\underline A$ associated with $A$ admits the following concrete description: for every $U\in X\etal$, the set $\underline A(U)$ is the set of locally constant functions from $U$ to $A$; if $U$ is locally connected ({\em e.g} $U$ is an algebraic variety over a field) then it can also be described as $A^{\pi_0(U)}$. Il $A$ is a an abelian group, then $\underline A$ is a sheaf of abelian groups. As usual, we will most of the time write $A$ instead of $\underline A$, and call it the {\em constant sheaf $A$}. 

\medskip
\noindent
{\bf The sheaf of global functions}. The assignment $U\mapsto \sch O_U(U)$ is a sheaf on $X\etal$ (we already knew that it was a sheaf for the Zariski topology). This means that for every \'etale $U\in X\etal$, and for every \'etale covering $(U_i\to U)$, any family $(s_i)$ with $s_i\in \sch O_{U_i}(U_i)$ for all $i$ and $s_{i|U_i\cap U_j}=s_{j|U_i\cap U_j}$ for all
$(i,j)$ arises from a uniquely determined $s\in \sch O_U(U)$.

\medskip
Let us see what it does mean in an interesting particular case. Suppose that $X=U=\spec k$ for some field $k$ and that $(U_i\to U)$ is a {\em single} map $\spec L\to \spec k$ for some finite Galois extension $L$ of $k$ with Galois group $G$; let $p_1$ and $p_2$ be the two projections $$\spec L\times_{\spec k}\spec L\to \spec L.$$The claim says that every $l\in L$ such that $p_1^*(l)=p_2^*(l)$ belongs to $k$.

Now one has an isomorphism between $L\otimes_kL$ and the ring of $G$-indexed finite sequences of elements of $L$,  which sends $(l\otimes \lambda)$ to $(lg(\lambda))_{g\in G}$; modulo this isomorphism, $p_1^*(l)=(l)_g$ and $p_2^*(l)=(g(l))_g$ for every $l\in L$. 

\medskip
Therefore one has $p_1^*(l)=p_2^*(l)$ if and only if $g(l)=l$ for every $g\in G$. Hence, our claim simply rephrases one of the fundamental results of Galois theory, that is, the equality $k=L^G$. 

\medskip
\noindent
{\bf Some comments.} We thus see that the theory of \'etale sheaves helps thinking to the equality $k=L^G$ in a quite geometrical way. More precisely, it follows from above that the fact, for an element of $L$, to be G-invariant, can be thought as a 'coincidence' condition; and the fact that such an element comes from $k$ can be thought as the result of a glueing process this coincidence condition has made possible. 

\medskip
This geometric vision of Galois phenomena has turned out to be very fruitful, especially as far as the so-called {\em descent theory} is concerned. The typical questions addressed by this theory is: {\em does a given object defined on $L$ and equipped with a Galois action come from $k$ ?} It is very helpful to think to it as a glueing problem (of objects, no more of functions), the existence of the Galois action corresponding precisely to the datum of 'a compatible system of isomorphisms on the intersections'. 

\medskip
\noindent
{\bf Coherent sheaves.} More generally, let $\sch F$ be a coherent (Zariski) sheaf on $X$. Remind that by the very definition of $X\etal$ every object $U$ of $X\etal$ comes equipped with an \'etale morphism $\pi: U\to X$; let us assign to such an $U$ the $\sch O_U(U)$-module $\pi^*\sch F(U)$. One can prove that this is a sheaf on $X\etal$, which we denote by $\sch F\etal$; note that the sheaf we have considered above is then nothing but $(\sch O_X)\etal$.

\medskip
One can prove that for every coherent sheaf $\sch F$ on $X$ and every integer $i$, the groups $\H^i(X_{\rm Zar}, \sch F)$ and $\H^i(X\etal, \sch F\etal)$ are canonically isomorphic. Hence as far as cohomology of coherent sheaves is concerned, \'etale cohomology doesn't contain more information than Zariski cohomology. 

\medskip
{\em Remark.} One often writes $\sch F$ instead of $\sch F\etal$, when the context is sufficiently clear to indicate that one doesn't work with the Zariski sheaf $\sch F$ but with the \'etale sheaf naturally associated with it. We will sometimes make such an abuse -- and write $\sch O_X$ for the sheaf $U\mapsto \sch O_U(U)$ on $X\etal$, as an example; the latter is also usually sometimes denoted by ${\Bbb G}_a$. 

\medskip
\noindent
{\bf The multiplicative group}. The assignment $U\mapsto \sch O_U(U)^*$ is a sheaf of abelian groups on $X\etal$; it is a subsheaf (of multiplicative monoids) of $\sch O_X$ which is often denoted by $\gm$. 

The 90'th theorem by Hilbert says that there is a natural isomorphism $\H^1(X_{\rm Zar}, \gm)\simeq \H^1(X\etal,\gm)$; this means that any '\'etale  line bundle' - that is, an \'etale sheaf of $\sch O_X$-modules locally isomorphic to $\sch O_X$ for the \'etale topology on $X$-- is already locally isomorphic to $\sch O_X$ for the {\em Zariski} topology. 

More generally, for every $n$ an \'etale sheaf of $\sch O_X$-modules which is locally isomorphic to $\sch O_X^n$ for the \'etale topology on $X$ is already locally isomorphic to $\sch O_X^n$ for the {\em Zariski} topology\footnote{This can also be expressed as a comparison between \'etale and Zariski $\H^1$, but with coefficients $\mathsf{GL}_n$ -- which is not an {\em abelian} group; we won't speak more about it here.}. When $X=\spec k$ for $k$ a field, the latter claim has a very concrete meaning: if $L$ is a finite Galois extension of $k$ and if $G$ is the Galois group of $L$, then every finite dimensional $L$-vector space endowed with an action of $G$ is $G$-isomorphic to $L\otimes_k V$ for some $k$-vector space $V$ (in other words, it admits a basis over $L$ whose vectors are invariant under $G$). This is an example of the aforementioned descent theorems. 

\medskip
\noindent
{\bf The sheaf of roots of unity}. Let $n$ be an integer. Let $\mu_n$ be the presheaf on $X\etal$ defined by the assignment $U\mapsto \{f\in \gm(U), f^n=1\}$; it easily seen to be a subsheaf of $\gm$; we denote it by $\mu_n$ (the 'sheaf of $n$-th root of unity'). 

\medskip
\noindent
{\em Remark.} If the ambient ground scheme need to be precised, one will write $A_X, {\Bbb G}_{m,X}, \mathsf{GL}_{n,X}$ and $\mu_{n,X}$ instead of $A, \gm, \mathsf{GL}_n$ and $\mu_n$.  

\medskip
\noindent
{\bf An important example of locally constant sheaf.} Let us assume that the integer $n$ is invertible on $X$ and let $\Phi\in \ZZ[X]$ be the $n$-th cyclotomic polynomial. One checks that for every open affine subscheme $U=\spec A$ of $X$ the finite map $V:=\spec A[T]/\Phi\to U$ is an \'etale cover, and that the map $\overline m\mapsto \overline T^m$ induce an isomorphism of sheaves $(\ZZ/n\ZZ)_V\simeq \mu_{n,V}$; hence the restriction of $\mu_n$ to $V$ is constant, whence we deduce that $\mu_n$ {\em is locally constant for the \'etale topology on $X$.}

\medskip
\noindent
{\bf Locally constant sheaves and representations of $\pi_1$}. Let $A$ be a finite set (resp. a finite abelian group) and assume that $X$ is a connected and non-empty scheme; let $\overline x$ be a geometric point of $X$. There is exactly the same phenomenon as in usual topology: one can establish a natural dictionary between sheaves (resp. sheaves of abelian groups) on $X\etal$ which are locally isomorphic to the constant sheaf $A$ and discrete ({\em i.e.} with open stabilizers) actions of $\pi_1(X,\overline x)$ on $A$, up to conjugation by a permutation (resp. an automorphism) of $A$. 

Let us give a very simple example of this in the arithmetical context. Let $X=\spec k$ with $k$ a field, and let $\overline x$ be given by a separable closure $k^s$ of $k$ whose Galois group is denoted by $G$; fix an integer $n$ which is not zero in $k$, and let $\zeta$ be a primitive $n$-th root of $1$ in $k^s$. The datum of $\zeta$ induces an isomorphism between the groups $\ZZ/n\ZZ$ and $\mu_n(k^s)$, whence a discrete action of $\mathsf G$ on $\ZZ/n\ZZ$; this is precisely the action that corresponds to the sheaf $\mu_n$ on $\spec k \etal$, which has been seen to be locally isomorphic to the constant sheaf ${\ZZ/n\ZZ}$. 

\medskip
\noindent
{\bf Cohomology of a locally constant sheaf.} Let us return to a general, non-empty connected scheme $X$ endowed with a geometric point $\overline x$; let $A$ be a finite abelian group and let $\sch F$ be an \'etale sheaf on $X$ which is locally isomorphic to $A$; it corresponds to a discrete action of $\pi_1(X,\overline x)$ on $A$. 

\medskip
It can then be proven that $\H^1(X\etal,\sch F)\simeq \H^1(\pi_1(X,\overline x), A)$, where $A$ is endowed with the aforementioned $\pi_1(X,\overline x)$-action (see section \ref{coho} for the definition of $\H^1(\pi_1(X,\overline x),.)$). 

\medskip
{\em The particular case where $\sch F= A$.} In that case, the action of $\pi_1(X,\overline x)$ on $A$ is trivial, and $\H^1(\pi_1(X,\overline x), A)$ is nothing but $\mathsf{Hom}_{\rm cont}(\pi_1(X,\overline x),A)$ (the topology on $A$ being the discrete one). 

\medskip
The group $\mathsf{Hom}_{\rm cont}(\pi_1(X,\overline x),A)$ classifies up to isomorphism the couples $(Y\to X,\iota)$ where $Y\to X$ is a connected Galois cover and $\iota$ an embedding $\mathsf{Aut}\;(Y/X)\hookrightarrow A$; that is therefore also the case for $\H^1(X\etal,A)$.

\medskip
\noindent
{\bf A notation.} If $X$ is a scheme and if $\ell$ is a prime number, we will denote by $\shl X$ the category of sheaves on $X\etal$ which are locally isomorphic to a constant sheaf associated with a finite abelian $\ell$-group. 

\medskip
\noindent
{\bf Pull-backs}. Let $f:Y\to X$ is a morphism of schemes. To every sheaf $\sch F$ on $X\etal$ can be associated in a natural way a sheaf $f^*\sch F$ on $Y\etal$;  If $\sch F=A_X$ for some set $A$ then $f^*\sch F=A_Y$; one has $f^*\mu_{n,X}=\mu_{n,Y}$ for every $n$ invertible on $X$. 

If $\sch F$ is locally isomorphic to a constant sheaf $A_X$, then $f^*\sch F$ is locally isomorphic to $A_Y$ (note in particular that if  $\sch F\in \shl X$ for some prime number $\ell$, then $f^*\sch F\in \shl Y$). If one assumes moreover that $A$ is finite and that $Y$ and $X$ are non-empty and connected, one has a nice description of $f^*\sch F$, as follows. Let $\overline y$ be a geometric point of $Y$ and let $\overline x$ denote the corresponding geometric point of $X$. The locally constant sheaf $\sch F$ is given by an action of $\pi_1(X,\overline x)$ on $A$; its pull-back $f^*\sch F$ is then the locally constant sheaf on $Y\etal$ that corresponds to the action of $\pi_1(Y,\overline y)$ on $A$ obtained by composing that of $\pi_1(X,\overline x)$ with $\pi_1(Y,\overline y)\to \pi_1(X,\overline x)$.

\medskip
If $X$ is a scheme over a field $k$, if $L$ is an extension of $k$ and if $\sch F$ is an \'etale sheaf on $X$, its pull-back to $X_L$ will be denoted by $\sch F_L$

\medskip
Assume now that $X$ is a scheme of finite type over $\CC$ and let $\sch F$ be an \'etale sheaf on $X$. There is a natural (topological) sheaf $\sch F_{\rm top}$ on $X(\CC)$ associated to $\sch F$. If $\sch F=A_X$ for some set $A$ then $\sch F_{\rm top}=A_{X(\CC)}$.

 If $\sch F$ is locally isomorphic to a constant sheaf $A_X$, then $\sch F_{\rm top}$ is locally isomorphic to $A_{X(\CC)}$; if one assumes moreover that $A$ is finite and that $X$ is connected and non-empty, there is a nice description of $\sch F_{\rm top}$ as follows. Let $x\in X(\CC)$ and let us still denote by $x$ the corresponding geometric point of $X$.  The locally constant sheaf $\sch F$ is given by an an action of $\pi_1(X,x)$ on $A$. the sheaf $\sch F_{\rm top}$ is then the locally constant sheaf that corresponds to the action of $\pi_1(X(\CC),x)$ on $A$ obtained by composing that of $\pi_1(X, x)$ with the arrow $\pi_1(X(\CC),x)\to\pi_1(X, x)$ (remind that the right-hand term is the profinite completion of the left-hand one). 

\subsection*{Two fundamental exact sequences} 

\noindent
{\bf The Kummer exact sequence.} Let $X$ be a scheme and let $n$ be an integer which is invertible on $X$. The {\em Kummer sequence} of \'etale sheaves on $X$ is the sequence $$ 1\to \mu_n\to \gm \stackrel{z\mapsto z^n}\longrightarrow \gm\to 1.$$ This is an exact sequence. To see it, the only point to check is that every section of $\gm$ has locally for the \'etale topology a preimage by $z\mapsto z^n$, that is, a $n$-th square root. 

\medskip
So, let $U\in X\etal$ and let $f$ an invertible function on $U$. We would like to prove that there exists an \'etale covering $(U_i)$ of $U$ such that $f_{|U_i}$ is an $n$-th power for every $i$; one immediately reduces to the case where $U$ is affine, say equal to $\spec A$. But $\spec A[T]/(T^n-f)\to U$ is then an \'etale covering (see the examples in section \ref{etmorsch}), and $\overline T$ is a $n$t-th root of $f$ in  $A[T]/(T^n-f)$. 

\medskip
\noindent
{\em Remark.} The Kummer exact sequence is often used as an algebraic substitute for the exponential sequence. 

\medskip
\noindent
{\bf The Artin-Schreier exact sequence.} Let $X$ be a scheme over $\FF_p$. The map $\sch O_X\stackrel{z\mapsto z^p-z}\longrightarrow \sch O_X$ is $\FF_p$ linear, and it is not difficult to see that its kernel is precisely $\underline {\FF_p}$. The {\em Artin-Schreier sequence} of \'etale sheaves on $X$ is the sequence $$0\to \underline{\FF_p}\to \sch O_X\stackrel{z\mapsto z^p-z}\longrightarrow \sch O_X\to 0.$$ This is an exact sequence. To see it, it remains to check that every section of $\sch O_X$ has locally for the \'etale topology a preimage by $z\mapsto z^p-z$. 

\medskip
So, let $U\in X\etal$ and let $f$ a function on $U$. We would like to prove that there exists an \'etale covering $(U_i)$ of $U$ such that $f_{|U_i}$ is equal for every $i$ to $z_i^p-z_i$ for a suitable function $z_i$ on $U_i$ ; one immediately reduces to the case where $U$ is affine, say equal to $\spec A$. But $\spec A[T]/(T^p-T-f)\to U$ is then an \'etale covering (see the examples in section \ref{etmorsch}), and $f=\overline T^p-\overline T$ in  $A[T]/(T^p-T-f)$. 

\subsection*{\'Etale cohomology with support}

We would like to mention here two important notions related to the \'etale cohomology of schemes.

\medskip
\noindent
{\bf The cohomology with support}. Let $X$ be a scheme and let $Y$ be a Zariski-closed subset of $X$. The functor that sends a sheaf of abelian groups $\sch F$ on $X\etal$ to the group of sections of $\sch F$ with support in $Y$ (that is, which restrict to zero on $X-Y$) is left-exact, hence admits derived functors. The latter are denoted by $\H^i_Y(X\etal,.)$; we say that  the $\H^i_Y(X\etal,\sch F)$'s are the {\em \'etale cohomology groups with support in $Y$} of $\sch F$. 

\medskip
\noindent
{\bf The cohomology with proper (or compact) support}. Let $X$ be a scheme, and assume that {\em $X$ is separated and of finite type over a field $k$}. One can then define for every sheaf $\sch F$ on $X\etal$ the {\em \'etale cohomology groups with proper (or compact) support of $\sch F$}, which are denoted by $\H^\bullet_c(X\etal,\sch F)$. 

The group $\H^0_c(X\etal,\sch F)$ is precisely the group of sections of $\sch F$ with proper support (that is, whose restriction to $X-Z$ is zero for some suitable closed subset $Z$ of $X$ which is proper), and if $X$ is itself proper then $\H^i_c(X\etal,\sch F)=\H^i(X\etal,\sch F)$ for every $i$ and every $\sch F$, as expected. 

\medskip
But we would like to emphasize that in general, contrary to the definition that is used classical in topology, the $\H^i_c(X\etal,.)$'s are {\em not} the derived functors of $\H^o_c(X\etal,.)$, though they give rise to long exact sequences starting from short ones. 

Those derived functors do exist, but are in general simply irrelevant, because there are 'not enough' closed subsets of $X$ which are proper  -- if $X$ is affine every such closed subset is a union of finitely many closed points, for example; that's why Grothendieck suggested a more {\em ad hoc} definition. 

\medskip
\noindent
{\small 
{\bf The definition of $\H^i_c(X\etal,\sch F)$}. As $X$ is of finite type and separated over a field $k$, it admits a {\em compactification}, that is, an open immersion $j$ into a proper $k$-scheme $\overline X$. Now define $j_{!}\sch F$ as the 'extension by zero' of $\sch F$, that is, the sheafification on $\overline X\etal$ of the presheaf that sends any $U\in \overline X\etal$ to $\sch F(U)$ if $U\to \overline X$ goes through $X$ and to zero otherwise.  

\medskip
One can prove that the groups $\H^i(\overline X\etal,j_{!}\sch F)$ don't depend on the chosen compactification of $X$, and for every $i$ one {\em sets} $$\H^i_c(X\etal,\sch F)= \H^i(\overline X\etal,j_{!}\sch F).$$

\medskip
{\em Remark.} One can prove that if $T$ is a compact topological space, if $\iota : U\hookrightarrow T$ is an open immersion, and if $\sch G$ is a sheaf on $U$ one has natural isomorphisms $$\H^i_c(U,\sch G)\simeq \H^i(T,\iota_{!}\sch G)$$ where $\iota_{!}$ denotes the 'extension by zero' functor; this gives kind of a motivation for the above definition. }

\medskip
\noindent
{\bf Link between \'etale and Zariski topology.} If $X$ is a scheme and if $\sch F$ is an \'etale sheaf of abelian groups on $X$, we know that $\sch F$ is also a sheaf for the Zariski topology on $X$. It follows from the very definition of derived functors (which we didn't give in full detail) that there are for every $i$ natural homomorphisms $\H^i(X_{\rm Zar},\sch F)\to \H^i(X\etal,\sch F)$; the same also holds for cohomology with support in a given Zariski-closed subset of $X$.  

\subsection*{\'Etale cohomology of an algebraic variety over an algebraically closed field}

\medskip
\noindent
{\bf Theorem 3}. {\em Let $k$ be an algebraically closed field and let $X$ be a separated $k$-scheme of finite type. Let $\ell$ be a prime number, and let $\sch F\in \shl X$. 

\medskip
A) If $\dim {}X=d$ then $\H^i(X\etal,\sch F)=0$ for $i>2d$, and if $X$ is affine then $\H^i(X\etal,\sch F)=0$ for $i>d$. 

\medskip
B) For every $i$ the group $\H_c^i(X\etal,\sch F)$ is finite, and if $L$ is an algebraically closed extension of $k$ then $\H_c^i(X\etal,\sch F)\to \H_c^i(X_{L,{\rm \acute e t}}, \sch F_L)$ is an isomorphism. 

\medskip
C) Assume that $X$ is proper or that $\ell\neq 0$ in $k$. Under those assumptions $\H^i(X\etal,\sch F)$ is finite, and if $L$ is an algebraically closed extension of $k$ then $\H^i(X\etal,\sch F)\to \H ^i(X_{L,{\rm \acute e t}}, \sch F_L)$ is an isomorphism. 

\medskip
D) If $k=\CC$ one has for every $i$ natural isomorphisms $$\H^i(X\etal,\sch F)\simeq \H^i(X(\CC),\sch F_{\rm top})\;{\rm and}\;\H_c^i(X\etal,\sch F)\simeq \H_c^i(X(\CC),\sch F_{\rm top}).$$}

\medskip
\noindent
{\em Remark.} There is a complex-analytic analog of assertion A): if $X$ is a complex analytic space of dimension $d$, then $\H^i(X,\ZZ)=0$ for $i>2d$ (because the topological dimension of $X$ is $2d$), and if $X$ is Stein, then  $\H^i(X,\ZZ)=0$ for $i>d$ (because one can prove that $X$ has then the homotopy type of a $d$-dimensional CW-complex); it is a general fact that affine schemes play a role analog to that of  Stein spaces. 

\medskip
\noindent
{\bf Good behaviour.} Assertion D) ensures that if one works with a scheme of finite type over $\CC$, then the \'etale cohomology groups will be exactly those one expects. But this also remains true over an arbitrarily algebraically closed field, at least for some broad classes of varieties. 

This claim seems probably quite vague to the reader. Let us be more precise through three examples, that of a curve, that of a torus and that of an abelian variety. We keep the notations $k$ and $\ell$ and we assume that $\ell\neq 0$ in $k$. We fix an integer $n$ and we set $\Lambda=\ZZ/\ell^n\ZZ$ 

\medskip
$\bullet$ If $X$ is a smooth, projective, irreducible curve of genus $g$, then $\H^1(X\etal,\Lambda)$ is isomorphic to $\Lambda^{2g}$. If $n>0$, if $P_1,\ldots,P_n$ are distinct $k$-points on $X$ and if $U=X-\{P_1,\ldots,P_n\}$, one has an exact sequence $$0\to \H^1(X\etal, \Lambda)\simeq \Lambda^{2g}\to\H^1(U\etal,\Lambda)\to \bigoplus_{i=1}^n \Lambda \stackrel \Sigma\longrightarrow \Lambda \to 0.$$

\medskip
$\bullet$ Let $n\in \NN$ and set $X={\Bbb G}_{m,k}^n$. The group $ \H^1(X\etal,\Lambda)$ is isomorphic to $\Lambda^n$; the connected Galois covers of $X$ which are associated to elements of  $ \H^1(X\etal,\Lambda)$ are exactly the isogenies $X'\to X$ (with connected $X'$) whose kernel embeds in $\Lambda$. 

The algebra $\H^\bullet(X\etal,\Lambda)$ (multiplication is the cup-product) is isomorphic to the free exterior $\Lambda$-algebra generated by $ \H^1(X\etal,\Lambda)\simeq \Lambda ^n$. 

\medskip
$\bullet$ Let $X$ be an abelian variety of dimension $g$ over $k$.  The group $ \H^1(X\etal,\Lambda)$ is isomorphic to $\Lambda^{2g}$;  the connected Galois covers of $X$ which are associated to elements of  $ \H^1(X\etal,\Lambda)$ are exactly the isogenies $X'\to X$ (with connected $X'$) whose kernel embeds in $\Lambda$. 

The algebra $\H^\bullet(X\etal,\Lambda)$ is isomorphic to the free exterior $\Lambda$-algebra generated by $ \H^1(X\etal,\Lambda)\simeq \Lambda ^{2g}$. 

\medskip
\noindent
{\bf Bad behaviour.} We will now show that  C) definitely fails if $\ell=0$ in $k$ and if $X$ is non proper. For that purpose, assume that char. $k=p>0$ and that $\ell=p$ (for the moment, we don't make any particular assumption on $X$). On $X\etal$ the Artin-Schreier exact sequence $$0\to \FF_p\to \sch O_X\stackrel{z\mapsto z^p-z}\longrightarrow \sch O_X\to 0$$ induces a long exact sequence $$0\to \H^0(X\etal,\FF_p)\to \H^0(X\etal,\sch O_X)\to \H^0(X\etal,\sch O_X)\to \H^1(X\etal,\FF_p)\to \H^1(X\etal,\sch O_X)\to \ldots$$

\medskip
\noindent
{\em Remark.} Assume that $X$ is connected and non-empty. Let $f\in \H^0(X\etal, \sch O_X)$; to its image in $\H^1(X\etal,\FF_p)$ is associated a connected Galois cover $Y$ of $X$, together with an embedding of $\mathsf{Aut}\;Y/X)$ into $\FF_p$; let us describe concretely $Y$ and the embedding. 

\medskip
- if $f=z^p-z$ for some $z\in  \H^0(X\etal,\sch O_X)$ then $Y=X$ (it is the case where the image of $f$ in $\H^1(X\etal,\FF_p)$ is zero), and the embedding is of course the trivial one; 

- if not, let $Y_U$ denote, for $U=\spec A$ an open affine subset of $X$, the scheme $\spec A[T]/T^p-T-f$; this is an \'etale cover of $U$, and the required Galois cover $Y\to X$ is obtained by glueing the $Y_U$'s for varying $U$. The embedding is then an isomorphism $\mathsf{Aut}\;Y/X\simeq \FF_p$, defined as follows: if $\phi$ is an $X$-automorphism of $Y$, there exists a (unique) $\lambda_\phi\in \FF_p$  such that $\phi^*\overline T=\overline T+\lambda_\phi$; the assignment $\phi\mapsto\lambda_\phi$ is the required isomorphism. 

\medskip
Suppose now that $X=\Aff^1_k$. By the comparison theorem for coherent sheaves, on has $\H^1(X\etal,\sch O_X)\simeq \H^1(X_{\rm Zar}, \sch O_X)$; as $X$ is affine, the latter is zero. We thus have, in view of the above long exact sequence, a natural isomorphism $$\H^1(X\etal,\FF_p)\simeq k[T]/\{z^p-z,z\in k[T]\}.$$ We thus see that $\H^1(X\etal,\FF_p)$ is infinite, and that it actually increases by extension of the algebraically closed ground field; hence both assertions of C) are false in this setting.

\medskip
\noindent
{\small
{\bf A few words about the proof of assertion D)}. Let us begin with a general remark. Let $\mathsf C$ be a site, let $\sch G$ be a sheaf on $\mathsf C$ and let $U$ be an object of $\mathsf C$; fix a cohomolgy class $h\in \H^\bullet (U, \sch G)$. One can then prove that, exactly like what happens in the classical case, the class $h$ is locally trivial on $U$ (for the given Grothendieck topology of $\mathsf C$, of course). 

Now let us go back to assertion D). It implies, in view of the aforementioned fact, the following assertion:

\medskip
(*)~~~~~ if $h\in H^\bullet(X(\CC),\sch F_{\rm top})$, then there exists an \'etale covering $(X_i\to X)$ such that the pull-back of $h$ to $X_i(\CC)$ is zero for all $i$. 

\medskip
In fact, by kind of a very general 'abstract non-sense', it turns out that to prove assertion D) (at least for cohomology without support), it is {\em sufficient} to prove (*). That is the way one does it when $X$ is smooth, by using Riemann's existence theorem; the general case is proven with Hironaka's resolution of singularities. }

\medskip
\subsection*{The $\ell$-adic formalism}

Let us keep the notations $k,X$ and $\ell$ of th. 3, and assume that $\ell\neq 0$ in $k$. Though they are perfectly well-defined, the groups $\H^i(X\etal,\ZZ_\ell)$ and $\H^i(X\etal,\QQ_\ell)$ in the sense of derived functors are irrelevant; let us explain why through a basic example. 

If $X$ is connected an non-empty then  $\H^1(X\etal,\QQ_\ell)\simeq \mathsf{Hom}_{\rm cont}(\pi_1(X),\QQ_\ell)$; as $\pi_1(X)$ is profinite and as $\QQ_\ell$ has no torsion, the latter is zero; therefore  $\H^1(X\etal,\QQ_\ell)=0$, which is not reasonable ($X$ is here {\em any} algebraic $k$-variety). 

\medskip
But as Grothendieck pointed out, the groups $$\lim_{\stackrel \leftarrow n} \H^i(X\etal,\ZZ/\ell^n\ZZ)\;\:{\rm and}\;\; \left( \lim_{\stackrel \leftarrow n}\H^i(X\etal,\ZZ/\ell^n\ZZ)\right)\otimes_{\ZZ_\ell}\QQ_\ell$$ are relevant. {\bf For that reason, the notation  $\H^i(X\etal,\ZZ_\ell)$ (resp. $\H^i(X\etal,\QQ_\ell)$) will almost always abusively denote the first one (resp. the second one). }

The same holds for cohomology with compact support. 

\medskip
With those abuses of notations theorem 3 above and the three examples which follow it also hold with $\sch F=\Lambda=\ZZ_\ell$ or $\QQ_{\ell}$, the notion of a finite group having to be replaced by that of a finitely generated $\ZZ_\ell$-module or $\QQ_\ell$-vector space (for assertion C) of th. 3 note that we suppose now $\ell\neq 0$ in $k$).  

\medskip
\noindent
{\bf An interesting consequence for $\ell$-adic coefficients. } Assume that $k=\CC$ and let $\sigma$ be any (that is, possibly non-continuous) automorphism of $\CC$. The corresponding twisted variety $X_\sigma$ coincides with $X$ as a scheme; only its structure map to $\spec \CC$ has been changed (through $\sigma$). 

Applying D) to $X$ and $X_\sigma$ and using the fact that they coincide as schemes yields to the existence of natural isomorphism of graded algebras $$\H^\bullet(X(\CC),\ZZ_\ell)\simeq \H^\bullet(X_\sigma(\CC),\ZZ_\ell)$$  and $$ \H^\bullet(X(\CC),\QQ_\ell)\simeq \H^\bullet(X_\sigma(\CC),\QQ_\ell).$$

\medskip
Note that the second isomorphism immediately implies that $X(\CC)$ and $X_\sigma(\CC)$ have the same Betti numbers. But when $X$ is smooth and projective\footnote{If $X$ is only assumed to be smooth, this is still possible using Grothendieck's GAGA results for de Rham cohomology, which involve resolution of singularities.}, one can prove it directly without using \'etale cohomology: Serre observed that it is a consequence of GAGA and of the fact that the Betti numbers of a complex projective (or, more generally, compact K\"ahler) manifold $Y$ can be computed, through Hodge theory, from the dimensions of the $\H^q(Y,\Omega^p)$'s. 

\medskip
\noindent
{\bf A counter-example for real coefficients}. Recently, Fran\c cois Charles has given in \cite{fch} an example of a complex, smooth projective variety $X$ and of an automorphism $\sigma$ of $\CC$ such that the graded {\em algebras} $\H^\bullet(X(\CC),\RR)$ and $\H^\bullet(X_\sigma(\CC),\RR)$ are not isomorphic; but note that the graded {\em vector spaces} $\H^\bullet(X(\CC),\RR)$ and $\H^\bullet(X_\sigma(\CC),\RR)$ are isomorphic, because $X(\CC)$ and $X_\sigma(\CC)$ have the same Betti numbers: the problem comes from the cup-product. 

\section{\'Etales sheaves and \'etale  cohomology of analytic spaces}

{\bf \em General references: {\rm \cite{brk2}} and {\rm \cite{brkc}} (the latter is used only for the comparison theorem -- that is, th. 5). }

{\em We fix a non-Archimedean, complete field $k$.}

\subsection*{Examples of \'etale sheaves}

\medskip
Let us begin with a remark: as open immersions are \'etale, any open covering of an analytic space $X$ is an \'etale covering of $X$; therefore every sheaf on $X\etal$ is, when restricted to the category of open subsets of $X$, a usual sheaf. 

 \medskip
We will quickly mention some important \'etale sheaves; as we will see, there is a great similarity between these examples and those which were considered in the scheme-theoretic context. 

\medskip
\noindent
{\bf Constant sheaves}. If $A$ is any set, the sheaf on $X\etal$ which is associated to the presheaf $U\mapsto A$ sends every $U\in X\etal$ to the set of locally constant $A$-valued functions on $U$ -- which coincides with $A^{\pi_0(X)}$; as usual, it will be most of the time simply denoted by $A$, when there is no risk of confusion. 

\medskip
As for schemes, for every prime number $\ell$ we will denote by $\shl X$ the category of \'etale sheaves on $X$ that are locally isomorphic to a constant sheaf associated with a finite abelian $\ell$-group.

\medskip
\noindent
{\bf Coherent sheaves}. Assume for that paragraph that $X$ is {\em good} (i.e. every point of $X$ has an affinoid neighborhood) and let $\sch F$ be a coherent sheaf of $\sch O_X$-modules on $X$. The assignment that sends every $\pi:U\to X$ in $X\etal$ to $\pi^*\sch F(U)$ is a sheaf on $X\etal$, which is denoted by $\sch F\etal$ or by $\sch F$ if this doesn't cause any trouble. For that kind of sheaves, \'etale cohomology doesn't carry more information than the topological one: one has for every $i$ a canonical isomorphism $\H^i(X\etal,\sch F)\simeq \H^i(X_{\rm top},\sch F)$. 

\medskip
\noindent
{\small {\bf What happens in the non-good case ?} What we have just said above remains essentially true, up to the following: coherent sheaves and their cohomology have to be considered with respect not to the usual topology of $X$ (it is too coarse) but to the G-topology (see the definition at section \ref{grotop}) ; as we already mentioned, this is a Grothendieck topology, but far more set-theoretic than the \'etale one. 

In the good case, the theories of coherent sheaves (including their cohomology groups) for the usual topology on $X$ and for the admissible one turn out to be the same.}

\medskip
\noindent
{\bf The multiplicative group}. The assignment $U\mapsto \sch O_U(U)^*$ is a sheaf on $X\etal$, which is usually denoted by $\gm$. If $X$ is good, one has a natural isomorphism $\H^1(X\etal,\gm)\simeq \H^1(X_{\rm top},\gm)$; in the general case, as above, one has to replace the topological $\H^1$ with the $\H^1$ related to the G-topology on $X$. 

\medskip
\noindent
{\bf The sheaf of roots of unity.} For every $n\in \NN$ the assignment $$U\mapsto \{f\in \sch O_U(U)^*, f^n=1\}$$ is a subsheaf of $\gm$ which is denoted by $\mu_n$. If $n\neq 0$ in $k$, this sheaf is locally isomorphic to the constant sheaf $\ZZ/n\ZZ$, and the {\em Kummer sequence} $$ 1\to \mu_n\to \gm \stackrel{z\mapsto z^n}\longrightarrow \gm\to 1$$ is exact. 

\medskip
\noindent
{\bf Locally constant sheaves.} Let $A$ be a finite set (resp. a finite abelian group) and assume that $X$ is connected and non-empty ; let $\overline x$ be a geometric point of $X$. As for schemes, one can establish a natural dictionary between sheaves (resp. sheaves of abelian groups) on $X\etal$ which are locally isomorphic to the constant sheaf $A$ and discrete ({\em i.e.} with open stabilizers) actions of $\pi_{1,\rm f}(X,\overline x)$ on $A$, up to conjugation by a permutation (resp. an automorphism) of $A$. 

\medskip
If $A$ is an abelian group and if $\sch F$ is the locally constant sheaf associated with a given action of $\pi_{1,\rm f}(X,\overline x)$ on $A$ then $\H^1(X\etal,\sch F)\simeq \H^1(\pi_{1,\rm f}(X,\overline x), A)$; in the particular case where $\sch F=A$ this is nothing but $\mathsf{Hom}_{\rm cont}(\pi_{1,\rm f}(X,\overline x),A)$, which classifies up to isomorphism the couples $(Y\to X,\iota)$ where $Y\to X$ is a connected Galois cover and $\iota$ an embedding $\mathsf{Aut}\;(Y/X)\hookrightarrow A$.

\medskip
\noindent
{\bf Pull-backs }. If $f:Y\to X$ is a morphism of analytic spaces, to every sheaf $\sch F$ on $X\etal$ can be associated in a natural way a sheaf $f^*\sch F$ on $Y\etal$. If $\sch F=A_X$ for some set $A$ then $f^*\sch F=A_Y$; one has $f^*\mu_{n,X}=\mu_{n,Y}$ for every $n$ invertible on $X$. 

If $\sch F$ is locally isomorphic to a constant sheaf $A_X$, then $f^*\sch F$ is locally isomorphic to $A_Y$ (note in particular that if $\sch F\in \shl X$ for some prime number $\ell$, then $f^*\sch F\in \shl Y$). If  one assumes moreover that $A$ is finite and that $Y$ and $X$ are non-empty and connected, one has a nice description of $f^*\sch F$, as follows. Let $\overline y$ be a geometric point of $Y$ and let $\overline x$ denote the corresponding geometric point of $X$. The locally constant sheaf $\sch F$ is given by an action of $\pi_{1,\rm f}(X,\overline x)$ on $A$; its pull-back $f^*\sch F$ is then the locally constant sheaf on $Y\etal$ that corresponds to the action of $\pi_{1,\rm f}(Y,\overline y)$ on $A$ obtained by composing that of $\pi_{1,\rm f}(X,\overline x)$ with $\pi_{1,\rm f}(Y,\overline y)\to \pi_{1,\rm f}(X,\overline x)$.

\medskip
If $X$ is a $k$-analytic space, if $L$ is a complete extension of $k$ and if $\sch F$ is an \'etale sheaf on $X$, its pull-back to $X_L$ will be denoted by $\sch F_L$. 

\medskip
If $\sch X$ is a $k$-scheme of finite type, then to every sheaf $\sch F$ on $\sch X\etal$ is associated in natural way a sheaf $\sch F\an$ on $\sch X\an \etal$. If $\sch F=A_{\sch X}$ for some set $A$ one has then $\sch F\an=A_{\sch X\an}$; if $\sch F=\mu_{n,\sch X}$ for some $n$ invertible on $\sch X$ then $\sch F\an=\mu_{n,\sch X\an}$.

If $\sch F$ is locally isomorphic to a constant sheaf $A_{\sch X}$, then $\sch F\an$ is locally isomorphic to $A_{\sch X}\an$ (note in particular that if $\sch F\in \shl {\sch X}$ for some prime number $\ell$, then $\sch F\an\in \shl {\sch X\an}$). If  one assumes moreover that $A$ is finite and that $\sch X$, and hence $\sch X\an$, are non-empty and connected, one has a nice description of $\sch F\an$, as follows. Let $\overline x$ be a geometric point of $\sch X\an$ and let $\overline \xi$ denote the corresponding geometric point of $\sch X$. The locally constant sheaf $\sch F$ is given by an action of $\pi_1(\sch X,\overline \xi)$ on $A$; its pull-back $\sch F\an$ is then the locally constant sheaf on $\sch X\an\etal$ that corresponds to the action of $\pi_{1,\rm f}(\sch X\an,\overline x)$ on $A$ obtained by composing that of $\pi_1(\sch X,\overline \xi)$ with $\pi_{1,\rm f}(\sch X\an,\overline x)\to \pi_1(\sch X,\overline \xi)$.

\subsection*{\'Etale cohomology with support}

As for schemes, one can define on $X$ cohomology groups with support in a given closed subset (the definition will be the same), and also, if $X$ is Hausdorff, cohomology groups with compact support -- it will be different, and in fact simpler. 

\medskip
\noindent
{\bf Cohomology with support}. If $T$ is a closed subset of $X$, the functor that sends a sheaf of abelian groups $\sch F$ on $X\etal$ to the group of sections of $\sch F$ with support in $T$ (that is, which restrict to zero on $X-T$) is left-exact, hence admits derived functors. The latter are denoted by $\H^i_T(X\etal,.)$; we say that  the $\H^i_T(X\etal,\sch F)$'s are the {\em \'etale cohomology groups with support in $T$} of $\sch F$. 

\medskip
\noindent
{\bf Cohomology with compact support.} One can also define cohomology with compact support for \'etale sheaves on a Hausdorff analytic space. Contrary to what happens in scheme theory, it is done like in topology, that is, by deriving the functor 'global sections with compact support'. 

\medskip
\noindent
{\bf Link between \'etale and usual topology.} If $X$ is an analytic space and if $\sch F$ is an \'etale sheaf of abelian groups on $X$, we know that $\sch F$ is also a sheaf for the usual topology on $X$. It follows from the very definition of derived functors (which we didn't give in full detail) that there are for every $i$ natural homomorphisms $\H^i(X_{\rm top},\sch F)\to \H^i(X\etal,\sch F)$; the same also holds for cohomology with support in a given closed subset of $X$.

\subsection*{\'Etale cohomology of analytic spaces: two general theorems}

\medskip
\noindent
{\bf Theorem 4 (the case of an algebraically closed ground field)}. {\em Assume that $k$ is algebraically closed and let $\ell$ be a prime number which is non-zero in $\red k$. Let $X$ be a Hausdorff $k$-analytic space and let $\sch F\in \shl X$. 

1) If $X$ is paracompact of dimension $d$, then $\H^i(X\etal,\sch F)=0$ for $i>2d$. 

2) If $L$ is an algebraically closed complete extension of $k$, then for every $i$ the natural arrows $\H^i (X\etal, \sch F)\to \H^i (X_{L,{\rm \acute e t}}, \sch F_L)$ and $\H^i_c (X\etal, \sch F)\to \H^i_c (X_{L,{\rm \acute e t}}, \sch F_L)$ are isomorphisms. 

3) Assume that $X$ is a finite union of affinoid domains, each of which can be embedded in the analytification of an algebraic $k$-variety. The group $\H^i (X\etal, \sch F)$ is finite.}

\medskip
Let us make some comments about the rather technical assumption on $X$ which is made in assertion 3). 

\medskip
$\bullet$ It is automatically true as soon as the space $X$ is compact and quasi-smooth\footnote{This means essentially means that $X$ can be G-locally defined by a system of equations satisfying the usual jacobian condition; it may have a boundary ({\em e.g.} any closed polydisc is quasi-smooth); a space is called {\em smooth} if it is quasi-smooth and boundaryless ({\em e.g} an open polydisc).} ({\em rig-smooth} in the rigid-analytic language). 

$\bullet$ Compactness of $X$ is likely sufficient for  $\H^i (X\etal, \sch F)$ to be finite; but at the moment local algebraicity is technically needed for proving it -- the only known way to do it being the reduction to some algebraic finiteness results (due to Deligne) through Berkovich's comparison theorem on vanishing cycles.

\medskip
\noindent
{\bf Theorem 5 (comparison between algebraic and analytic \'etale cohomology groups).} {\em Let $\sch X$ be a separated $k$-scheme of finite type, let $\ell$ be a prime number and let $\sch F\in \shl {\sch X}$. 

1) One has $\H^i_c(\sch X\an\etal,\sch F\an)\simeq \H^i_c (\sch X\etal,\sch F)$.

2) If moreover the prime $\ell$ is non zero in $k$ (but possibly zero in $\red k$) then $\H^i(\sch X\an\etal,\sch F\an)\simeq \H^i (\sch X\etal,\sch F)$.}   

\medskip
\noindent
{\bf About $\ell$-adic formalism in analytic geometry.} If $\ell$ is a prime number which is non zero in $k$ then on a compact $k$-analytic space one can define relevant \'etale cohomology groups with coefficients $\ZZ_ \ell$ or $\QQ_\ell$ exactly in the same way than in scheme theory. For non-compact spaces it is slightly more complicated;  the only reference on that topic are, as far as we know, unpublished notes by Berkovich. 

\subsection*{The cohomology of some analytic spaces}

We will give some examples which show that the cohomology of analytic spaces over an algebraically closed field 'takes the expected values'; here it will mean that the values will be the same as they are for the analog complex analytic space. {\em For that section, we suppose that $k$ is algebraically closed}, we fix a prime number $\ell$ which is non zero in $\red k$, an integer $n$ and we set $\Lambda=\ZZ/\ell^n\ZZ$. 

\medskip
$\bullet$  If $X$ in an open polydisc, then $\H^i (X\etal,\Lambda)=0$ for every $i>0$. 

\medskip
$\bullet$ If $X$ is an open poly-annulus whose dimension is denoted by $d$, then $\H^1(X\etal,\Lambda)\simeq \Lambda^d$,; the connected Galois covers of $X$ which are associated to elements of  $ \H^1(X\etal,\Lambda)$ are exactly the connected Kummer covers $X'\to X$ whose group embeds in $\Lambda$. 

The algebra $\H^\bullet(X\etal,\Lambda)$ is isomorphic to the free exterior $\Lambda$-algebra generated by  $\H^1(X\etal,\Lambda)\simeq \Lambda^d$.

\medskip
$\bullet$ If $X$ is a connected, proper, smooth analytic group of dimension $g$ ({\em e.g.} the analytification of an abelian variety, or $\gm^{g,an}/\mathsf L$ for some lattice $\mathsf L\subset (k^*)^g$)  then $\H^1(X,\Lambda)\simeq \Lambda^{2g}$; the connected Galois covers of $X$ which are associated to elements of  $ \H^1(X\etal,\Lambda)$ are exactly the isogenies $X'\to X$ (with connected $X'$) whose kernel embeds in $\Lambda$. 

The algebra $\H^\bullet(X\etal,\Lambda)$ is isomorphic to the free exterior $\Lambda$-algebra generated by  $\H^1(X\etal,\Lambda)\simeq \Lambda^{2g}$.

\section{Poincar\'e duality}

{\bf \em General references: {\rm \cite{sga43}}, {\rm \cite{sga45}} and {\rm \cite{mil}} in the case of schemes; {\rm \cite{brk2}} for analytic spaces. }

\medskip
Poincar\'e duality theorems in both algebraic and analytic settings look formally the same; we will therefore use a common formulation. 

\noindent
{\bf Cup-products with supports.} Let us begin with a general remark, which hold for topological spaces as well as for  \'etale topology of schemes and analytic spaces. In those theories, cup-product pairings are also defined for cohomology with a given closed support, and for cohomology with compact support when the latter makes sense. They give rise to pairings $$\H^i_Y(.,.) \times \H^j_Z(.,.)\to \H^{i+j}_{Y\cap Z}(.,.)$$ and $$\H^i_c(.,.) \times \H^j (.,.)\to \H^{i+j}_c(.,.)$$ which behave as expected.

\medskip
\noindent
{\bf Poincar\'e duality: context and notations.} Now, let us denote by $k$ an algebraically closed field (resp. an algebraically closed non-Archimedean field). Let $X$ be a separated, smooth, $k$-scheme of finite type (resp. an Hausdorff smooth $k$-analytic space) of pure dimension $d$. We fix a prime number $\ell$ which is non zero in $k$ (resp. $\red k$) and an integer $n$; set $\Lambda=\ZZ/\ell^n\ZZ$. 

If $\sch F\in \shl X$ is annihilated by $\Lambda$ and if $i\in \NN$, we will denote by $\sch F(i)$ the sheaf $\sch F\otimes\mu_{\ell^n}^{\otimes i}$, and by $\sch F^\vee$ the sheaf that sends $U$ to ${\rm Hom}(\sch F_{|U}, \Lambda_U)$ (it is called the {\em dual} sheaf of $\sch F$); if $i\in \ZZ-\NN$ the sheaf $\sch F(-i)^\vee$ will be denoted by $\sch F(i)$. The formulas $\sch F(i+j)=\sch F(i)\otimes\sch F(j)$ and $\sch F(i)^\vee=\sch F(-i)$ then hold for all $i, j$ in $\ZZ$. 

\medskip
\noindent
{\bf Theorem 5 (Poincar\'e duality for schemes and analytic spaces).} {\em We use the notations above. 

1) There is a natural {\em trace} map $\mathsf{Tr}: \H^{2d}_c(X\etal,\Lambda(d))\to \Lambda$.

2)  For every $\sch F\in \shl X$ which is annihilated by $\Lambda$ and every $i\in\{0,\ldots,2d\}$ the map $$\H^i(X\etal,\sch F)\times\H^{2d-i}_c(X\etal, \sch F^\vee(d))\to \Lambda, $$ induced by the cup-product pairing and the trace map, identifies $\H^i(X\etal,\sch F)$ with the dual as $\Lambda$-module of  $\H^{2d-i}_c(X\etal, \sch F^\vee(d))$, that is, with $$\mathsf{Hom}_\Lambda(\H^{2d-i}_c(X\etal, \sch F^\vee(d)),\Lambda).$$}

\medskip
\noindent
{\bf Let us make some comments.} 

\medskip
{\em A particular fundamental case.} If $X$ is connected and non-empty, then $\H^0(X\etal,\Lambda)\simeq \Lambda$; we thus deduce of 2) that the trace map induces an isomorphism $\H^{2d}_c(X\etal,\Lambda(d))\simeq \Lambda$. 

\medskip
{\em About duality.} In the algebraic case, duality assertion 2) only involves {\em finite} $\Lambda$-modules. But that is not the case in the analytic context. Let us give a very simple example of that phenomenon in the {\em zero-dimensional} case (which is, maybe surprisingly, not completely empty): if $X$ is a disjoint union of one-point spaces $\sch M(k)$ indexed by some set $I$ (possibly infinite) then $\H^0(X,\Lambda)=\Lambda^I$ and $\H^0_c(X,\Lambda)=\Lambda^{(I)}$; classical duality between those two groups is then a particular case of Poincar\'e duality ! 

\medskip
{\em About comparison with  transcendent Poincar\'e duality.} In general topology, Poincar\'e  duality between the cohomology without support and that with compact support on a purely-dimensional topological variety involves the {\em orientation sheaf} of the variety (which is locally isomorphic to the constant sheaf $\ZZ$) ; in \'etale cohomology, the role of the latter is, in view of the theorem above, played by the sheaf of roots of unity. 

But while working with smooth complex analytic spaces, the orientation sheaf is usually non mentioned: it is because every smooth complex analytic space
$X$ comes with a natural orientation (defined through multiplication by $i$), whence an isomorphism between the orientation sheaf of $X$ and $\ZZ_X$. But we emphasize that this identification is related to the choice of one of the two square-roots of $-1$ in $\CC$; such a choice also induces for every $N$ an isomorphism $\mu_N\simeq \ZZ/N\ZZ$ (by distinguishing $e^{2i\pi/N}$), modulo which the isomorphisms of th. 3, assertion D) are compatible with algebraic and transcendent Poincar\'e dualities for $N$-torsion finite coefficients. 

\medskip
\noindent
{\bf About the definition of the trace map.} We keep the notations of the theorem above.  Let $D$ be an irreducible and reduced divisor on $X$, that is, a reduced closed subscheme (resp. analytic subspace) of dimension $d-1$. As $X$ is smooth, $D$ is Cartier ({\em i.e.} locally defined by one non-zero equation). The sheaf of rational (or meromorphic) functions defined on $X-D$ and having at most simple poles along $D$ is then a line bundle $\sch O(D)$ with a canonical section $1$, which generates $\sch O(D)$ on $X-D$. To such a datum is associated in a quite formal way a cohomology class of $\H^1_D(X_{\rm top}, \gm)$ which one can push to $\H^1_D(X\etal,\gm)$,  and then to $\H^2_D(X\etal,\mu_{\ell^n})=\H^2_D(X\etal,\Lambda(1))$ through the long-exact cohomology sequence associated with the Kummer exact sequence of \'etale sheaves $$1\to \mu_{\ell^n}\to \gm \stackrel{z\mapsto z^{\ell^n}}\longrightarrow \gm \to 1.$$

This construction extends to an assignment which sends every Zariski closed subset (resp. every closed analytic subset) $Y$ of $X$ of pure codimension $i$ to a cohomology class in $\H_Y^{2i}(X\etal,\Lambda(i))$, which can be pushed to $\H^{2i}_c(X\etal,\Lambda(i))$ if $Y$ is proper. The formation of those cohomology classes associated with Zariski-closed (resp. closed analytic) subsets of $X$ makes {\em transverse} intersections commute with cup-products. 

\medskip
Now let $x$ be a $k$-point on $X$. By the above, one can associate to it a class $h(x)\in \H^{2d}_c(X\etal,\Lambda(d))$. One then has $\mathsf{Tr} (h(x))=1\in \Lambda$.  As a consequence, if $X$ is non-empty and connected, then $h(x)$ doesn't depend on $x$ and generates the free $\Lambda$-module $ \H^{2d}_c(X\etal,\Lambda(d))$; it is called the {\em fundamental class} of $X$.

\section{The case of a non-algebraically closed ground field}

{\bf \em General references: {\rm \cite{sga42}} or {\rm \cite{mil}} for algebraic varieties, and {\rm \cite{brk2}} for analytic spaces. }

Up to now most of our theorems about the cohomology of algebraic varieties or analytic spaces concerned the case of an {\em algebraically closed} field. We would like to give some explanations about what happens if one removes that assumption. 

\subsection*{\'Etale cohomology of a field and Galois cohomology}

Let $k$ be a field, let $k^s$ be a separable closure of $k$ and let $G$ be the profinite group $\mathsf{Gal}\;(k^s/k)$. If $\sch F$ is an \'etale sheaf on $\spec k$ and if $L$ is a finite sub-extension of $k^s/k$, then $\spec L\to \spec k$ is \'etale; therefore $\sch F(\spec L)$ does make sense, and will be written $\sch F(L)$ for the sake of simplicity. The direct limit of the $\sch F(L)$'s, for $L$ going through the set of finite sub-extensions of $k^s/k$, is in a natural way a discrete $G$-set (discrete means that the stabilizers are open); we denote it by $\sch F(k^s)$; if $\sch F$ is a sheaf in abelian groups, then $\sch F(k^s)$ is a discrete $G$-module. 

\medskip
\noindent
{\bf Theorem 6.} {\em We use the above notations.

1) The assignment $\sch F\mapsto \sch F(k^s)$ induces an {\em equivalence} (i.e. a categorical 1-1 correspondence) between the category of \'etale sheaves on $\spec k$ and that of discrete $G$-sets. 

2) If $\sch F$ is an \'etale sheaf on $\spec k$ and if $L$ is a finite sub-extension of $k^s/k$ corresponding to an open subgroup $H$ of $G$ then $\sch F(L)=\sch F(k^s)^H$. 

3) If $\sch F$ is a sheaf of abelian groups on $\spec k$ then one has for every $i$ an isomorphism $\H^i((\spec k)\etal,\sch F)\simeq \H^i(G,\sch F(k^s))$. }

\medskip
\noindent
{\bf Some comments.} We therefore see that the \'etale cohomology theory of $\spec k$ is nothing but the cohomology theory of discrete $\mathsf{Gal}\;(k^s/k)$-modules, which is also called the {\em Galois cohomology} of the field $k$. Galois cohomology plays a crucial role in modern arithmetic -- we suggest the interested reader  to look at Serre's book (\cite{cgl}).

\medskip
\noindent
{\bf A topological analog}. Let $X$ be a connected, non-empty, 'reasonable' topological space, admitting a universal cover $Y\to X$; let $G$ be the automorphism group of $Y$ over $X$. Let $A$ be an abelian group equipped with an action of $G$. It corresponds to a local system $\sch F$ on $X$ which is locally isomorphic to $A$; the pull-back of $\sch F$ on $Y$ is {\em globally} isomorphic to $A$. The group $G$ acts in a natural way on $\H^\bullet(Y,A)$. 

\medskip
There is a natural spectral sequence $$\H^p(G,\H^q(Y,A))\Rightarrow \H^{p+q}(X,\sch F).$$ Now assume that $\H^q(Y,A)=0$ for every $q>0$ ({\em e.g} $Y$ is contractible). Then this spectral sequences simply furnishes isomorphisms $\H^p(X,\sch F)\simeq \H^p(G,A)$ for every $p$. Therefore we see that, as far as \'etale cohomology is concerned, the spectrum of a field $k$ behaves, to some extent, like a connected, non-empty topological space having a contractible universal cover; the spectrum of a separable closure $k^s$ of $k$ can be thought as such a cover. 

\subsection*{Varieties or analytic spaces over a field}

Our purpose is now to explain how the cohomology of a variety or an analytic space over a field $k$ mixes, in some sense, Galois cohomology of the ground field and cohomology of the given variety over an algebraic closure of $k$. 

Let us begin with some remarks. Fix a field $k$ and an algebraic closure $k^a$ of $k$ ; let $G$ be the absolute Galois group $\mathsf{Gal}\;(k^a/k)$.

1) If $X$ is an algebraic variety over $k$ and if $\sch F$ is an an \'etale sheaf of abelian groups on $X$, there is a natural discrete action of $G$ on $\H^\bullet(X_{k^a,{\rm \acute e t}},\sch F_{k^a})$. 

2) If $k$ is non-Archimedean and complete, and if $\widehat{k^a}$ denotes the completion of $k^a$, then for every {\em compact}\footnote{The fact for the action of $G$ on the cohomology of $\sch F_{\widehat{k^a}}$ to be discrete essentially means that every cohomology class of  $\sch F_{\widehat{k^a}}$ is already defined on a finite separable extension of $k$; one needs compactness to ensure this finiteness result.} analytic space $X$ and every \'etale sheaf of abelian groups $\sch F$ on $X$, there is a natural discrete action of $G$ on $\H^\bullet(X_{\widehat{k^a},{\rm \acute e t}},\sch F_{\widehat{k^a}})$. 

\medskip
Now we will see that in the first (resp. second) case, the cohomology of $\sch F$ is builded from the $\H^p(G,\H^q(X _{k^a,{\rm \acute e t}},\sch F_{k^a}))$'s (resp.  $\H^p(G,\H^q (X _{\widehat{k^a},{\rm \acute e t}},\sch F_{\widehat{k^a}})$'s) through a spectral sequence. 

\medskip
Let us now state a precise theorem.

\medskip
\noindent
{\bf Theorem 7.} {\em Let $k$ be a field and let $k^a$ be an algebraic closure of $k$; let $G$ be the absolute Galois group $\mathsf{Gal}\;(k^a/k)$. 

\medskip
1) If $X$ is a $k$-scheme of finite type and if $\sch F$ is an \'etale sheaf of abelian groups on $X$, one has a natural spectral sequence $$\H^p (G,\H^q(X_{k^a, {\rm \acute e t}},\sch F_{k^a}))\Rightarrow \H^{p+q}(X\etal,\sch F).$$

\medskip
2) Assume that $k$ is non-Archimedean and complete, and let $\widehat{k^a}$ be the completion of $k^a$. If $X$ is a {\em compact} analytic space over $k$ and if $\sch F$ is an \'etale sheaf of abelian groups on $X$, one has a natural spectral sequence $$\H^p (G,\H^q(X_{\widehat{k^a}, {\rm \acute e t}},\sch F_{\widehat{k^a}}))\Rightarrow \H^{p+q}(X\etal,\sch F).$$}

 \end{document}